\documentclass[11pt,reqno]{amsart}
\usepackage{amsmath, amsfonts, amsthm, mathtools, extarrows,amssymb,mathrsfs,dsfont}
\usepackage[x11names]{xcolor}
\usepackage[colorlinks=true,linkcolor=blue,citecolor=blue,urlcolor=black]{hyperref}
\usepackage{cleveref}
\usepackage{tikz-cd}
\usepackage[shortlabels]{enumitem}
\setlist[enumerate]{leftmargin=*}
\usepackage[margin=1.05in]{geometry}
\usepackage[indent,skip=3pt plus 1pt minus 1pt]{parskip}
\usepackage[style=alphabetic,backref,url=false,maxnames=10,giveninits,isbn=false]{biblatex}
\newtheorem{theorem}{Theorem}[section]
\newtheorem{lemma}[theorem]{Lemma}
\newtheorem{corollary}[theorem]{Corollary}
\newtheorem{proposition}[theorem]{Proposition}

\theoremstyle{definition}
\newtheorem{definition}[theorem]{Definition}
\newtheorem{remark}[theorem]{Remark}

\numberwithin{equation}{section}

\usepackage[T1]{fontenc}

\setitemize[1]{itemsep=5pt}
\setenumerate[1]{itemsep=5pt}

\DeclarePairedDelimiterX\set[1]\lbrace\rbrace{#1}

\title[Injectivity of the Kudla--Millson lift]
{The injectivity of the Kudla--Millson lift on unitary Shimura varieties}
\author{Hanyu Yao}
\address{Shanghai Center for Mathematical Sciences, Fudan University, Shanghai 200438, China}
\email{hyyao23@m.fudan.edu.cn}
\date{}
\subjclass[2020]{Primary 11F27; Secondary 11F37, 11F55, 14G35}
\keywords{Borcherds products, Kudla--Millson lift, unitary Shimura varieties, harmonic weak Maass forms, Heegner divisors}
\begin{document}
\begin{abstract}
We prove the injectivity of the Kudla--Millson lift on unitary Shimura varieties over imaginary quadratic fields under some lattice hypotheses.
We express integrals of the lift over special cycles as special values of Dirichlet series and deduce injectivity by separating the Fourier coefficients of the cusp forms.
We also establish a converse theorem for Borcherds products on $\mathrm{U}(n,1)$.
\end{abstract}
\maketitle

\setcounter{tocdepth}{1}
\tableofcontents
\section{Introduction}

The Kudla--Millson lift \cite{KudlaMillsonHarmonicI,KudlaMillsonHarmonicII,KudlaMillsonTubes} connects cusp forms with the geometry of special cycles on orthogonal and unitary arithmetic quotients.
The associated theta series define holomorphic modular forms in cohomology, with Fourier coefficients given by the Poincar\'e dual classes of Heegner cycles \cite{KM90}.
In this paper, we study the injectivity of the Kudla--Millson lift on ball quotients, viewed as connected components of unitary Shimura varieties over imaginary quadratic fields.
Under a Hermitian splitting hypothesis, we show that the action of units on the discriminant module accounts for the entire kernel of the lift.

Another motivation for this injectivity problem comes from the theory of Borcherds products.
Borcherds's multiplicative lift \cite{borcherds1996automorphic} constructs meromorphic automorphic forms on orthogonal Shimura varieties from weakly holomorphic modular forms with integral principal part.
Their divisors are linear combinations of Heegner divisors determined by the Fourier coefficients of the weakly holomorphic forms.
The converse problem asks when a meromorphic automorphic form whose divisor is a linear combination of Heegner divisors arises from this construction.
In the orthogonal setting, the work of Bruinier \cite{Bruinier2002book} and Bruinier--Funke \cite{BFtwolift,BruinierFunkeKMInjectivity} relates this question to the injectivity of the Kudla--Millson lift.

In the unitary setting, Hofmann \cite{hofmann2014borcherds} constructed Borcherds products by restricting orthogonal products along the trace embedding. 
Meanwhile, Funke and Hofmann \cite{FunkeHofmannGreen} developed the underlying singular theta lift and Green current theory, and Hofmann \cite{hofmann2022FJ} later computed the corresponding Fourier--Jacobi expansion.

\subsection{Main results}
Let $L$ be a Hermitian lattice in a Hermitian space $(V,h)$ of signature $(n,1)$ over an imaginary quadratic field $K$.
We assume that the trace lattice $L_{\mathbb Z}$ is even.
Let $A_L$ be its discriminant module and let $\rho_L$ be the associated finite Weil representation of $\operatorname{SL}_2(\mathbb Z)$.
The discriminant kernel $\Gamma_L^{\mathrm U}\subseteq\operatorname{U}(L)$ acts on the complex $n$-ball $\mathcal D_{V}$ of negative lines in $V\otimes_K\mathbb C$, yielding the quotient $X_L=\Gamma_L^{\mathrm U}\backslash\mathcal D_{V}$ and its finite covers $X_\Gamma=\Gamma\backslash\mathcal D_{V}$ for finite-index subgroups $\Gamma\subseteq\Gamma_L^{\mathrm U}$.

The Kudla--Millson theta function defines a conjugate linear lift $\Lambda_{\mathrm{U}}:S_{n+1,\rho_L}\to\mathcal Z^{1,1}(X_{\Gamma})$ into the space of smooth closed $(1,1)$-forms on $X_{\Gamma}$. 
The natural action of $\mathcal U=\mathcal O_K^\times$ on $A_L$ induces an action on $S_{n+1,\rho_L}$.
Let $P_{\mathcal U}$ denote the averaging projection onto the $\mathcal U$-invariant subspace $S_{n+1,\rho_L}^{\mathcal U}$.

Take $\delta_K\in\mathcal O_K$ such that $\delta_K^2=\operatorname{disc}(K)$.
In this paper, we say that $L$ splits an orthogonal dual hyperbolic plane over $\mathcal{O}_K$ if $L=L_0\perp H_{\mathrm{tr}}$, where $H_{\mathrm{tr}}=\mathcal O_Ke\oplus\mathcal O_Ke'$ satisfies $h(e,e)=h(e',e')=0$ and $h(e,e')=\delta_K^{-1}$.
When $n\geq2$, this condition is equivalent to requiring that $L_{\mathbb Z}$ split an $\mathcal O_K$-stable orthogonal direct summand isometric to $2U$, where $U$ is the even unimodular hyperbolic plane over $\mathbb Z$.
\begin{theorem}
\label{thm:introKMKernel}
Assume that $n>2$ and that $L$ splits an orthogonal dual hyperbolic plane over $\mathcal{O}_K$.
Then $\ker(\Lambda_{\mathrm{U}})=\ker(P_{\mathcal U})$.
In particular, $\Lambda_{\mathrm{U}}$ is injective on $S_{n+1,\rho_L}$ if and only if $P_{\mathcal U}=\operatorname{id}$.
\end{theorem}

Under the hypotheses of Theorem \ref{thm:introKMKernel}, $\Lambda_{\mathrm U}$ is injective whenever $K\notin\{\mathbb Q(i),\mathbb Q(\sqrt{-3})\}$.
Subsection \ref{subsec:kmKernel} establishes more refined injectivity criteria for $K=\mathbb Q(i)$ and $\mathbb Q(\sqrt{-3})$, while Subsection \ref{subsec:eisensteinLatticeExamples} provides geometric examples satisfying these hypotheses. 
Moreover, the counterexample in Subsection \ref{subsec:splittingCounterexample} shows that our splitting hypothesis cannot in general be replaced by an abstract $2U$-splitting of $L_{\mathbb Z}$ or by the existence of an $\mathcal O_K$-stable orthogonal summand isometric to $U\perp U(N)$ with $N>1$.

\begin{remark}
Let $\rho_L^{\mathrm U}$ be the finite Weil representation of $\mathrm U(1,1)(\mathbb Z)$.
Then $S_{n+1,\rho_L^{\mathrm U}}=S_{n+1,\rho_L}^{\mathcal U}$, so the unit-invariance condition is part of the modular transformation law.
Since $\Lambda_{\mathrm U}=\Lambda_{\mathrm U}\circ P_{\mathcal U}$, Theorem~\ref{thm:introKMKernel} equivalently states that $\Lambda_{\mathrm U}$ is injective on $S_{n+1,\rho_L^{\mathrm U}}$.
We retain the $\operatorname{SL}_2(\mathbb Z)$ formulation to match the modular form spaces used in the references on unitary Borcherds products and in our converse argument.
\end{remark}
We also prove a weak converse theorem for unitary Borcherds products without any splitting hypothesis.
\begin{theorem}
\label{thm:introWeakConverse}
Assume that $n\geq2$ and $F$ is a meromorphic modular form for $\Gamma$ with a finite order multiplier system, whose divisor is a linear combination of Heegner divisors.
Then there exists a harmonic weak Maass form $f$  with associated Heegner divisor $\operatorname{div}(F)$ such that $\log\lVert F\rVert+\frac14\Phi_L^{\mathrm U}(\,\cdot\,,f)$ is constant on $X_\Gamma\setminus|\operatorname{div}(F)|$.
\end{theorem}
Here $\lVert F\rVert$ is the Petersson norm of $F$, and $\Phi_L^{\mathrm U}$ is the singular theta lift, which expresses the logarithmic Petersson norm of the Borcherds product.
The Heegner divisor associated with $f$ is the linear combination of Heegner divisors determined by the negative index Fourier coefficients of its holomorphic part, as defined in Subsection \ref{subsec:singularThetaLift}.

Theorem \ref{thm:introWeakConverse} serves as the unitary analogue of Bruinier's weak converse theorem \cite[Theorem 4.23]{Bruinier2002book}. 
Combined with our injectivity result, it yields a full converse theorem for unitary Borcherds products.
For simplicity, we state the result here only for fields with $\mathcal O_K^\times=\{\pm1\}$.
\begin{theorem}
\label{thm:introBorcherdsConverse}
For $n>2$, assume that $K\notin\{\mathbb Q(\sqrt{-1}),\mathbb Q(\sqrt{-3})\}$ and $L$ splits an orthogonal dual hyperbolic plane over $\mathcal O_K$.
Then every meromorphic modular form $F$ with a finite order multiplier system whose divisor is a linear combination of Heegner divisors is (up to a nonzero constant factor) the Borcherds product of a  weakly holomorphic modular form. 
\end{theorem}
\subsection{Methods of proof}
Several methods have been developed to establish the injectivity of orthogonal Kudla--Millson lifts under suitable lattice hypotheses. 
Bruinier initially proved injectivity by analyzing the Fourier expansion of the associated regularized theta lift \cite{Bruinier2002book}, later combining this strategy with a newform theory for the Weil representation \cite{Bruinier14Converse}. 
An alternative approach computes the Fourier expansion of the Kudla--Millson lift directly via unfolding. 
Introduced by Zuffetti \cite{ZuffettiKMUnfolding}, this technique was later extended to general signatures and Funke--Millson lifts by Metzler and Zuffetti \cite{MetzlerZuffettiKM}, and to the adelic setting over totally real fields by Li, Zhang, and Zuffetti \cite{LiZhangZuffettiTotallyRealKM}.

Taking a different route, Bruinier--Funke \cite{BruinierFunkeKMInjectivity} and Stein \cite{stein2023conversetheoremborcherdsproducts} established injectivity by applying the Rallis inner product formula to evaluate the $L^2$-norm of the lift. 
Most recently, Metzler \cite{metzler2025newconversetheoremborcherds, metzler2025symmetricsquaretypelseries} introduced a cycle integral method.
We adapt this strategy to the unitary setting.

To prove Theorem \ref{thm:introKMKernel}, we integrate $\Lambda_{\mathrm{U}}(g)$ for $g\in S_{n+1,\rho_L}$ over the special cycles attached to positive $K$-lines $V_1\subset V$. 
By decomposing the theta kernel along $V=V_1\perp V_1^\perp$ and applying the unitary Siegel--Weil formula together with Rankin--Selberg unfolding, we express these cycle integrals as special values of a Dirichlet series.

Recovering the unit averages of the Fourier coefficients of $g$ from these special values requires an additional arithmetic construction, as the rank-one projective $\mathcal O_K$-module $L^\vee\cap V_1$ need not be principal, unlike its orthogonal counterpart over $\mathbb Q$.
For each Fourier index $(m,\mu)$ of $g$ and each nonzero integral ideal $I\subseteq\mathcal O_K$, we choose a suitable pair of generators of $I^{-1}$.
Using these generators in the $H_{\mathrm{tr}}$-summand, we construct a positive line whose cycle integral determines the value at $s=n$ of the corresponding Dirichlet subseries summed over $I$.
If $\Lambda_{\mathrm{U}}(g)=0$, all these special values vanish.
Ideal M\"obius inversion then shows that the subseries over elements coprime to any nonzero integral ideal $J$ also vanishes at $s=n$.
Choosing $J$ divisible by all principal ideals $(\alpha)$ with $1<N(\alpha)\le R$ excludes these nonunit terms while retaining every unit term.
As $R\to\infty$, absolute convergence makes the remaining nonunit contribution tend to zero.
Thus the Fourier coefficients sum to zero over every $\mathcal U$-orbit, yielding $P_{\mathcal U}g=0$.

The proof of Theorem \ref{thm:introWeakConverse} involves an analytic feature distinct from the orthogonal case, as the complex ball is not a tube domain for $n\geq2$. 
While Bruinier, Howard, and Yang \cite{BruinierHowardYang2015} established the required $L^p$-integrability for the regularized theta lift, they stated that its Laplacian is constant away from the divisor by appealing to the orthogonal case without proof. 
We supply an explicit computation of this Laplacian. 
To control the growth of $F$ near the cusps, we use meromorphic continuation and toroidal coordinates to show that $\log\lVert F\rVert\in L^p(X_\Gamma)$ for $1<p<n$, which implies that the smooth function $\log\lVert F\rVert+\frac14\Phi_L^{\mathrm U}(\,\cdot\,,f)$ is constant.

Following the argument of \cite[Theorem 4.2]{Bruinier14Converse}, we combine the injectivity of the Kudla--Millson lift with Theorem \ref{thm:introWeakConverse} to obtain Theorem \ref{thm:introBorcherdsConverse}.

\subsection{Application to the Picard group}
These results also contribute to the study of relations among Heegner divisors in rational Picard groups.
For orthogonal Shimura varieties, the results of \cite{BMM17orthogonal} in the compact case and \cite{BLMM17} in the noncompact case show that special cycles generate rational Hodge classes in sufficiently low codimension.
In codimension one, these results imply that special divisors generate the rational Picard group.
Bruinier's work \cite{Bruinier2002PicardRank,Bruinier2002book,Bruinier14Converse} relates the dimension of the subspace of the rational Picard group spanned by Heegner divisors to the injectivity of the Kudla--Millson lift.
Since special divisors and Heegner divisors span the same subspace by \cite{li2026specialcyclespicardgroups}, the resulting dimension formulas give the full Picard rank.

In the unitary setting, the generation theorem of Bergeron, Millson, and Moeglin \cite{BMM16unitary} applies only to compact quotients and includes arbitrary rational $(1,1)$-classes among its generators, so it does not imply that special divisors generate the rational Picard group even in the compact case. 
Nevertheless, our results determine the dimension of the subspace $\operatorname{Pic}_{\mathrm{Heeg}}(X_L)_{\mathbb Q} \subset \operatorname{Pic}(X_L)\otimes_{\mathbb Z}\mathbb Q$ spanned by Heegner divisors.

\begin{corollary}
\label{cor:introHeegnerPicardDimension	}
Suppose that $n>2$ and that $L$ splits an orthogonal dual hyperbolic plane over $\mathcal{O}_K$.
Then
$$
\dim_{\mathbb Q}\operatorname{Pic}_{\mathrm{Heeg}}(X_L)_{\mathbb Q}=1+\dim_{\mathbb C}S_{n+1,\rho_L}^{\mathcal U}.
$$
In particular, if $P_{\mathcal U}=\operatorname{id}$ on $S_{n+1,\rho_L}$, this dimension equals $1+\dim_{\mathbb C}S_{n+1,\rho_L}$.
\end{corollary}

A further question concerns the Picard number of the Baily--Borel compactification.
In the orthogonal setting, Huang--Li--M\"uller--Ye \cite{HuangLiMullerYe2025PicardBailyBorel} prove that this number is one under suitable lattice and level hypotheses by using the local Borcherds products of Bruinier and Freitag \cite{BruinierFreitag2001Local} to study the extension of Heegner divisor classes across the boundary.
Since Hofmann \cite{hofmann2019local} constructed analogous local Borcherds products for unitary groups, one may ask whether the rational Picard classes on the unitary Baily--Borel compactification whose restrictions lie in $\operatorname{Pic}_{\mathrm{Heeg}}(X_L)_{\mathbb Q}$ likewise form a one-dimensional space under suitable lattice hypotheses.

\noindent\textbf{Outline.}
In Sections \ref{sec:unitary shimura varieties} and \ref{sec:weilModularForms}, we review the necessary background on unitary Shimura varieties, Heegner divisors and the Weil representations.
In Section \ref{sec:borcherdsConverse}, we recall unitary Borcherds products and Kudla--Millson lifts.
 We also state the weak converse theorem, which is proved in Section \ref{sec:weakConverseProof}, and use it to reduce the converse problem to the injectivity of the Kudla--Millson lift.
Section \ref{sec:cycleIntegrals	} evaluates the special cycle integrals of $\Lambda_{\mathrm{U}}$ in terms of Dirichlet series.
Finally, in Section \ref{sec:kmInjectivity}, we determine $\ker(\Lambda_{\mathrm{U}})$ and present our geometric examples and counterexample.

\noindent\textbf{Acknowledgements.}
The author is grateful to his advisor, Zhiyuan Li, for his constant support and many valuable comments on this work.
He also thanks Zelin Ye for helpful discussions.
The author is supported by NSFC Grant No.~12121001.

\noindent\textbf{History and AI disclosure.} 
This project originated in October 2025 from a suggestion of Zhiyuan Li to study the Picard number of the Baily--Borel compactification of unitary ball quotients. 
Through discussions with Li and Zelin Ye, it became clear that establishing the injectivity of the unitary Kudla--Millson lift was the first key step. 
Before July 2026, the author had completed a draft proof of the weak converse theorem following \cite{Bruinier2002book} and figured out how to overcome the difficulties from nonprincipal ideals specific to the unitary case, which made it possible to prove injectivity on the full adelic quotient under the $H_{\mathrm{tr}}$-splitting condition via cycle integrals.
At that stage, the author did not know how to prove the corresponding injectivity statement on a single connected component.
 
The author subsequently sought to relax this lattice hypothesis, as requiring the trace lattice to split an $\mathcal O_K$-stable orthogonal summand isometric to two hyperbolic planes appeared restrictive and it was also natural to ask whether $H_{\mathrm{tr}}$ could be replaced by a standard unimodular Hermitian hyperbolic plane.
Before adopting the cycle integral method, the author had tried to prove injectivity by computing Fourier--Jacobi expansions, but could not recover the individual Fourier coefficients of the cusp form.
Later attempts to return to this approach with LLM assistance also failed to yield a proof, even under the $H_{\mathrm{tr}}$-splitting assumption.
After these prolonged attempts, the author prompted an LLM to search for counterexamples under weaker lattice assumptions, which led to the example in Subsection \ref{subsec:splittingCounterexample}.

Aside from this, LLMs assisted in passing from the full adelic quotient to a single connected component in Proposition \ref{prop:adelicComponentIntegrals} and helped locate \cite{BruinierHowardYang2015}, which simplified the author's original proof of Proposition \ref{prop:weakConverseRegularity}.
In addition, LLMs were used for language writing, literature searches, and refining the exposition of certain mathematical arguments.
The author takes full responsibility for the mathematical content of this paper.
\section{Unitary Shimura varieties and Heegner divisors}
\label{sec:unitary shimura varieties}
\subsection{Unitary Shimura varieties}
\label{subsec:unitary shimura varieties}

Let $K=\mathbb Q(\sqrt{D_K})$ be an imaginary quadratic field of discriminant $D_K<0$, and fix an embedding $\sigma:K\hookrightarrow\mathbb C$.
Denote by $\delta_K\in K$ the unique element satisfying $\delta_K^2=D_K$ and $\operatorname{Im}\sigma(\delta_K)>0$.
The ring of integers of $K$ is $\mathcal O_K=\mathbb Z[\zeta_K]$, where
$$
\zeta_K=
\begin{cases}
\delta_K/2,&D_K\equiv0\pmod4,\\
(1+\delta_K)/2,&D_K\equiv1\pmod4.
\end{cases}
$$

Let $(V,h)$ be a Hermitian space over $K$ of dimension $n+1$ and signature $(n,1)$, where the Hermitian form $h$ is linear in the first variable. 
Set $V_{\mathbb C}=V\otimes_{K,\sigma}\mathbb C$. 
Equipping the underlying rational vector space $V_{\mathbb Q}=\operatorname{Res}_{K/\mathbb Q}V$ with the quadratic form $Q(x)=\frac12\operatorname{Tr}_{K/\mathbb Q}h(x,x)$ and the bilinear form $B(x,y)=\operatorname{Tr}_{K/\mathbb Q}h(x,y)$ yields a quadratic space $(V_{\mathbb Q},Q)$ of signature $(2n,2)$.

A Hermitian lattice in $V$ is a full projective $\mathcal O_K$-submodule $L\subset V$.
Its underlying $\mathbb Z$-module, equipped with $Q$ and $B$, is called the orthogonal trace lattice and denoted by $L_{\mathbb Z}$. 
A Hermitian lattice $L$ is called even if its trace lattice $L_{\mathbb Z}$ is even.
Throughout this paper, Hermitian lattices are assumed to be even.
The Hermitian dual and orthogonal trace dual are defined by
$$
L_h^\vee=\{x\in V:h(x,L)\subseteq\mathcal O_K\}
\quad\text{and}\quad
L^\vee=\{x\in V:B(x,L)\subseteq\mathbb Z\},
$$
which satisfy $L^\vee=\delta_K^{-1}L_h^\vee$.
The discriminant group $A_L=L^\vee/L$ carries the quadratic form $q_L\colon A_L\to\mathbb Q/\mathbb Z$ defined by $q_L(x+L)=Q(x)+\mathbb Z$.

Write $H=\mathrm U(V)$ for the unitary group of $V$ and $\mathrm U(L)$ for the subgroup preserving $L$.
The associated Hermitian symmetric domain is the complex $n$-ball $\mathcal D_{V}=\{[z]\in\mathbb P(V_{\mathbb C}):h(z,z)<0\}$.
Let $\mathbb A$ and $\mathbb A_f$ denote the adeles and finite adeles of $\mathbb Q$.
For a compact open subgroup $K_f\subset H(\mathbb A_f)$, the associated unitary Shimura variety has complex points
$$
X_{K_f}=H(\mathbb Q)\backslash(\mathcal D_{V}\times H(\mathbb A_f))/K_f.
$$
By \cite[Theorem 5.1]{BorelAdeles}, the double coset set $C_{K_f}=H(\mathbb Q)\backslash H(\mathbb A_f)/K_f$ is finite of cardinality $r=|C_{K_f}|$. 
Choosing representatives $h_1,\ldots,h_r\in H(\mathbb A_f)$ for $C_{K_f}$ and setting $\Gamma_j=\{\gamma\in H(\mathbb Q):h_j^{-1}\gamma_fh_j\in K_f\}$, we obtain a decomposition of $X_{K_f}$ into connected components via the maps $[z]\mapsto[(z,h_j)]$:
\begin{equation}
X_{K_f}\simeq\coprod_{j=1}^r\Gamma_j\backslash\mathcal D_{V}.
\label{eq:shimuraComponents}
\end{equation}
Define the unitary discriminant kernel by $\Gamma_L^{\mathrm U}=\ker\bigl(\mathrm U(L)\longrightarrow\mathrm O(A_L,q_L)\bigr)$.
Its adelic counterpart is the compact open subgroup $K_L=\ker\bigl(\operatorname{Stab}_{H(\mathbb A_f)}(\widehat L)\longrightarrow\operatorname{Aut}(A_L)\bigr)$, where $\widehat L=L\otimes_{\mathbb Z}\widehat{\mathbb Z}$.
The equality $V\cap\widehat L=L$ implies $H(\mathbb Q)\cap K_L=\Gamma_L^{\mathrm U}$, identifying the connected component of $X_{K_L}$ represented by $1\in H(\mathbb A_f)$ with the ball quotient $X_L=\Gamma_L^{\mathrm U}\backslash\mathcal D_V$.
Throughout the paper, we work with $X_L$ and its finite covers $X_\Gamma=\Gamma\backslash\mathcal D_V$ associated with finite index subgroups $\Gamma\subseteq\Gamma_L^{\mathrm U}$.
These quotients are complex orbifolds of dimension $n$.

\subsection{Orthogonal and unitary Heegner divisors}
\label{subsec:heegnerDivisors}
Let $(M,Q_M)$ be an even lattice of signature $(r,2)$ with bilinear form $B_M$ and discriminant form $(A_M,q_M)$.
Denote by $\mathcal D_M$ a connected component of the Grassmannian of oriented negative definite 
two-planes in $M\otimes_{\mathbb Z}\mathbb R$.
Given $\mu\in A_M$ and $m\in\mathbb Q_{>0}$ with $m+\mathbb Z=q_M(\mu)$, define the orthogonal Heegner divisor on $\mathcal D_M$ by
$$
\widetilde H_{\mathrm{O}}(m,\mu)
=\sum_{\substack{x\in\mu+M\\Q_M(x)=m}}
\{Z\in\mathcal D_M:B_M(z,x)=0\text{ for all }z\in Z\}.
$$
Similarly, for ${\mu}\in A_L$ and $m\in\mathbb Q_{>0}$ with $m+\mathbb Z=q_L({\mu})$, define the unitary Heegner divisor on $\mathcal D_V$ by
$$
\widetilde H_{\mathrm{U}}(m,{\mu})
=\sum_{\substack{x\in{\mu}+L\\Q(x)=m}}
\{Z\in\mathcal D_V:h(z,x)=0\text{ for all }z\in Z\}.
$$
Both sums are locally finite.

Write $\mathrm O^+(M)$ for the subgroup of $\mathrm O(M)$ preserving $\mathcal D_M$.
Define the orthogonal discriminant kernel by $\Gamma_M^{\mathrm O}=\ker\bigl(\mathrm O^+(M)\longrightarrow\mathrm O(A_M,q_M)\bigr)$ and put $X_M=\Gamma_M^{\mathrm O}\backslash\mathcal D_M$.
For any finite index subgroup $\Gamma\subseteq\Gamma_L^{\mathrm U}$, the divisors above are invariant under $\Gamma_M^{\mathrm O}$ and $\Gamma$, respectively, and descend to $H_{\mathrm{O}}(m,\mu)$ on $X_M$ and $H_{\mathrm{U}}(m,{\mu})$ on $X_\Gamma$.
The quotient divisors are understood as $\mathbb Q$-Cartier divisors.

We now specialize to $(M,Q_M)=(L_{\mathbb Z},Q)$ and identify $(A_M,q_M)$ with $(A_L,q_L)$.
The trace pairing induces an inclusion $\mathrm U(V)\hookrightarrow\operatorname{SO}(V_{\mathbb Q},Q)$.
Identifying $V_{\mathbb Q}\otimes_{\mathbb Q}\mathbb R$ with $V_{\mathbb C}$ as real vector spaces via $\sigma$, we view each negative complex line as a real two-plane oriented by $(z,-iz)$.
These planes lie in a single connected component $\mathcal D_{L_{\mathbb Z}}$ of the Grassmannian of oriented negative two-planes, and we choose $\mathcal D_M=\mathcal D_{L_{\mathbb Z}}$.
The resulting map
$$
\iota:\mathcal D_V\longrightarrow\mathcal D_{L_{\mathbb Z}},
\qquad
[z]\longmapsto\operatorname{span}_{\mathbb R}\{z,-iz\},
$$
is a holomorphic $\mathrm U(V)(\mathbb R)$-equivariant embedding.
The natural inclusion $\mathrm U(L)\subseteq\mathrm O^+(L_{\mathbb Z})$ restricts to $\Gamma_L^{\mathrm U}\subseteq\Gamma_{L_{\mathbb Z}}^{\mathrm O}$.
Consequently, $\iota$ induces a holomorphic map $\bar\iota:X_\Gamma\longrightarrow X_M$.
By \cite[Lemma 6.2]{hofmann2014borcherds} and descent to the quotients, we obtain
$\iota^*\widetilde H_{\mathrm{O}}(m,{\mu})=\widetilde H_{\mathrm{U}}(m,{\mu})$ and $\bar\iota^*H_{\mathrm{O}}(m,{\mu})=H_{\mathrm{U}}(m,{\mu})$.
\section{Weil representations and vector-valued modular forms}
\label{sec:weilModularForms}

In this section, we fix our conventions for the Weil representation and recall the basic notions of vector-valued modular forms and harmonic weak Maass forms.
\subsection{Adelic and finite Weil representations}
\label{subsec:weilRepresentations}

Writing $e(t)=e^{2\pi it}$, we fix the standard additive character $\psi:\mathbb Q\backslash\mathbb A\to\mathbb C^\times$ with archimedean component $\psi_\infty(t)=e(t)$.
Since $\psi$ is trivial on $\mathbb Q$, its finite part satisfies $\psi_f(r)=e(-r)$ for $r\in\mathbb Q$.

Let $(V',h')$ be a nondegenerate Hermitian space over $K$ of dimension $d$, with $h'$ linear in the first variable.
On the underlying rational vector space $V'_{\mathbb Q}$, we set $Q_{V'}(x)=h'(x,x)$ and $B_{V'}(x,y)=\operatorname{Tr}_{K/\mathbb Q}h'(x,y)$.
Let $W_0=K^2$ be the split skew-Hermitian plane with Gram matrix $\left(\begin{smallmatrix}0&1\\-1&0\end{smallmatrix}\right)$. 
Setting $G=\operatorname{U}(W_0)$ and $H'=\operatorname{U}(V')$, the natural matrix inclusion $\jmath:\operatorname{SL}_2\hookrightarrow G$ identifies $\operatorname{SL}_2$ with $\operatorname{SU}(W_0)$.
For $a\in\mathbb A_K^\times$ and $b\in\mathbb A$, write
$$
m(a)=\begin{pmatrix}a&0\\0&\overline a^{-1}\end{pmatrix},
\qquad
n(b)=\begin{pmatrix}1&b\\0&1\end{pmatrix},
\qquad
w=\begin{pmatrix}0&1\\-1&0\end{pmatrix}.
$$

Let $\epsilon_{K/\mathbb Q}$ be the quadratic idele class character attached to $K/\mathbb Q$, and choose a unitary character $\chi':K^\times\backslash\mathbb A_K^\times\to\mathbb C^\times$ whose restriction to $\mathbb A^\times$ is $\epsilon_{K/\mathbb Q}^{\,d}$.
The choice of $\psi$ and $\chi'$ determines a Weil representation $\omega_{V'}$ of the group $G(\mathbb A)\times H'(\mathbb A)$.
In the Schr\"odinger model on the space $\mathcal S(V'_{\mathbb Q}(\mathbb A))=\mathcal S(V'(\mathbb A_K))$ of Schwartz--Bruhat functions, it acts by
$$
\begin{aligned}
(\omega_{V'}(m(a),1)\varphi)(x)
&=\chi'(a)|a|_{\mathbb A_K}^{d/2}\varphi(ax),\\
(\omega_{V'}(n(b),1)\varphi)(x)
&=\psi\bigl(bQ_{V'}(x)\bigr)\varphi(x),\\
(\omega_{V'}(w,1)\varphi)(x)
&=\int_{V'_{\mathbb Q}(\mathbb A)}
\varphi(y)\psi\bigl(B_{V'}(x,y)\bigr)\,dy,\\
(\omega_{V'}(1,\gamma)\varphi)(x)
&=\varphi(\gamma^{-1}x),
\qquad \gamma\in H'(\mathbb A).
\end{aligned}
$$
Here $|\cdot|_{\mathbb A_K}$ is the idele norm, with $|a|_\infty=a\overline a$, and $dy$ is the self-dual Haar measure with respect to $B_{V'}$ and $\psi$.
For $\tau=u+iv\in\mathbb H$, put $g_\tau=n(u)m(\sqrt v)\in\operatorname{SL}_2(\mathbb R)$, viewed as an adelic element with trivial finite component.

To recall the finite Weil representation, let $J$ be a nondegenerate even $\mathbb Z$-lattice of signature $(b^+,b^-)$ and set $\sigma_J=b^+-b^-$.
Its discriminant group $A_J=J^\vee/J$ carries the quadratic form $q_J$ with associated bilinear form $b_J(\mu,\nu)=q_J(\mu+\nu)-q_J(\mu)-q_J(\nu)$.
Let $\widetilde T=(T,1)$ and $\widetilde S=(S,\sqrt\tau)$ be the standard lifts of $T=n(1)$ and $S=w^{-1}$ to the metaplectic double cover $\operatorname{Mp}_2(\mathbb Z)$ of $\operatorname{SL}_2(\mathbb Z)$.
The finite Weil representation $\rho_J=\rho_{A_J,q_J}$ on $\mathbb C[A_J]$ is determined on the standard basis $\{\mathbf e_\mu:\mu\in A_J\}$ by
$$
\begin{aligned}
\rho_J(\widetilde T)\mathbf e_\mu
&=e\bigl(q_J(\mu)\bigr)\mathbf e_\mu,\\
\rho_J(\widetilde S)\mathbf e_\mu
&=\frac{e(-\sigma_J/8)}{\sqrt{|A_J|}}
\sum_{\nu\in A_J}e\bigl(-b_J(\mu,\nu)\bigr)\mathbf e_\nu.
\end{aligned}
$$
This representation is unitary for the standard Hermitian inner product of $\mathbb C[A_J]$ and has finite image.
When $J$ has even rank, $\rho_J$ factors through $\operatorname{SL}_2(\mathbb Z)$.

Replacing $q_J$ with $-q_J$ yields the complex conjugate representation $\overline{\rho_J}$.
We write $\rho^\vee$ for the contragredient representation on the complex linear dual.
By unitarity, the map $\mathbf e_\mu^\vee\mapsto\mathbf e_\mu$ identifies $\rho_J^\vee$ with $\overline{\rho_J}$, where $\{\mathbf e_\mu^\vee\}$ is the dual basis.

Let $J\subset V'$ be a Hermitian lattice whose trace lattice is even, and let $\rho_J$ be the finite Weil representation of this trace lattice.
Here $J^\vee$ is the trace dual and $q_J(x+J)=Q_{V'}(x)+\mathbb Z$.
Put $\widehat J=J\otimes_{\mathbb Z}\widehat{\mathbb Z}$ and $\phi_\mu^J=\mathbf1_{\mu+\widehat J}$ for $\mu\in A_J$.
These characteristic functions form a basis of the subspace $\mathcal S_J^{\mathrm{fin}}=\bigoplus_{\mu\in A_J}\mathbb C\phi_\mu^J\subseteq\mathcal S(V'_{\mathbb Q}(\mathbb A_f))$.

\begin{proposition}[{\cite[(2.7)--(2.10)]{BruinierYang}}]
\label{prop:finiteSchwartzModel}
For $\gamma\in\operatorname{SL}_2(\mathbb Z)$ with finite adelic component $\gamma_f$, the operator $\omega_{V',f}(\jmath(\gamma_f),1)$ preserves $\mathcal S_J^{\mathrm{fin}}$.
The isomorphism
$$
\mathcal S_J^{\mathrm{fin}}\xrightarrow{\ \sim\ }\mathbb C[A_J]^\vee,
\qquad
\phi_\mu^J\longmapsto\mathbf e_\mu^\vee
$$
identifies this action with $\rho_J^\vee$.
The dual isomorphism therefore identifies $(\mathcal S_J^{\mathrm{fin}})^\vee$ with $\mathbb C[A_J]$ carrying the representation $\rho_J$.
\end{proposition}
\subsection{Vector-valued modular forms}
\label{subsec:vectorValuedModularForms}

Let $(W,\rho)$ be a finite-dimensional unitary representation of $\operatorname{SL}_2(\mathbb Z)$ such that $\rho(T)$ has finite order.
Set $j(\gamma,\tau)=c\tau+d$ for $\gamma=\left(\begin{smallmatrix}a&b\\c&d\end{smallmatrix}\right)\in\operatorname{SL}_2(\mathbb R)$ and $\tau\in\mathbb H$. 
Given $k\in\mathbb Z$ and a function $f:\mathbb H\to W$, the weight $k$ slash operator for $\gamma\in\operatorname{SL}_2(\mathbb Z)$ is defined by
$$
\bigl(f\vert_{k,\rho}\gamma\bigr)(\tau)
=j(\gamma,\tau)^{-k}\rho(\gamma)^{-1}f(\gamma\tau).
$$

\begin{definition}
\label{def:vectorValuedModularForms}
A holomorphic function $f:\mathbb H\to W$ is a vector-valued modular form of weight $k$ and type $\rho$ if it satisfies $f\vert_{k,\rho}\gamma=f$ for every $\gamma\in\operatorname{SL}_2(\mathbb Z)$ and is holomorphic at the cusp $\infty$.
Write $M_{k,\rho}$ for the space of these forms and $S_{k,\rho}$ for the subspace of cusp forms vanishing at $\infty$. 
Allowing a pole at $\infty$ gives the space $M^!_{k,\rho}$ of weakly holomorphic modular forms.
\end{definition}

For $f\in M^!_{k,\rho_L}$, the transformation law under $T$ yields the Fourier expansion
$$
f(\tau)=\sum_{\mu\in A_L} \sum_{\substack{m\in\mathbb Q\\m+\mathbb Z=q_L(\mu)}}c_f(m,\mu)e(m\tau)\mathbf e_\mu.
$$
When $f\in M^!_{k,\overline{\rho_L}}$, the summation runs instead over $m+\mathbb Z=-q_L(\mu)$. 
In both cases, $c_f(m,\mu)\neq 0$ for only finitely many $m<0$. 
Moreover, $c_f(m,\mu)=0$ for all $m<0$ if $f\in M_{k,\rho_L}$, and also for $m=0$ if $f\in S_{k,\rho_L}$.
In particular, the Fourier expansion of a cusp form $g\in S_{n+1,\rho_L}$ takes the form
$$
g(\tau)=\sum_{\mu\in A_L} \sum_{\substack{m>0\\m+\mathbb Z=q_L(\mu)}} b_g(m,\mu)e(m\tau)\mathbf e_\mu.
$$

For any finite quadratic module $A$, write $\{x,\lambda\}=\lambda(x)$ for the bilinear evaluation pairing between $\mathbb C[A]$ and $\mathbb C[A]^\vee$.
We equip $\mathbb C[A]$ with the standard Hermitian inner product $\langle\sum_{\mu\in A}a_\mu\mathbf e_\mu,\sum_{\mu\in A}b_\mu\mathbf e_\mu\rangle_A=\sum_{\mu\in A}a_\mu\overline{b_\mu}$, which is linear in the first argument.
\subsubsection*{Modular forms for $\mathrm U(1,1)(\mathbb Z)$}
We also consider vector-valued modular forms for the integral unitary group $\mathrm U(1,1)(\mathbb Z)=G(\mathbb Q)\cap\operatorname{GL}_2(\mathcal O_K)$ associated with $W_0$.
Write $\mathcal U=\mathcal O_K^\times$, viewed as a group of scalar matrices.
In this matrix model, one has $\operatorname{SU}(1,1)(\mathbb Z)=\operatorname{SL}_2(\mathbb Z)$ and $\mathrm U(1,1)(\mathbb Z)=\mathcal U\operatorname{SL}_2(\mathbb Z)$.

Take $V'=V$ in Subsection~\ref{subsec:weilRepresentations} and choose the splitting character with $\chi'_\infty(z)=(z/|z|)^{-(n+1)}$.
To obtain this normalization, start with a splitting character $\chi'_0$ and write $\chi'_{0,\infty}(z)=(z/|z|)^a$, where $a\equiv n+1\pmod2$.
By~\cite[Lemma~2.3.1]{PatrikisTate}, there exists a unitary Hecke character $\eta$ of $K$ with $\eta_\infty(z)=(z/|z|)^{-(n+1+a)/2}$.
Then $\chi'=\chi'_0\eta/\eta^c$, where $\eta^c(x)=\eta(\overline x)$, has the required infinity type and the same restriction to $\mathbb A^\times$.

For $u\in\mathcal U$, the triviality of $\chi'$ on $K^\times$ gives $\chi'_f(u)=\chi'_\infty(u)^{-1}=u^{n+1}$.
Since $m(u)=uI_2$ and $u\widehat L=\widehat L$, the formula for the action of $m(a)$ in Subsection~\ref{subsec:weilRepresentations} gives
$$
\omega_{V,f}(m(u),1)\phi_\mu^L=u^{n+1}\phi_{u^{-1}\mu}^L.
$$
Together with Proposition~\ref{prop:finiteSchwartzModel} and the group decomposition above, this shows that $\mathcal S_L^{\mathrm{fin}}$ is preserved by $\omega_{V,f}(\gamma_f,1)$ for every $\gamma\in\mathrm U(1,1)(\mathbb Z)$.
Under the dual identification of Proposition~\ref{prop:finiteSchwartzModel}, let $\rho_L^{\mathrm U}$ be the contragredient of this action on $\mathcal S_L^{\mathrm{fin}}$.
Its restriction to $\operatorname{SL}_2(\mathbb Z)$ is $\rho_L$, and the preceding calculation gives
$$
\rho_L^{\mathrm U}(u\gamma)\mathbf e_\mu
=u^{-(n+1)}\rho_L(\gamma)\mathbf e_{u^{-1}\mu}
$$
for $u\in\mathcal U$ and $\gamma\in\operatorname{SL}_2(\mathbb Z)$.
In particular, $\rho_L^{\mathrm U}$ is unitary and has finite image.
The centrality of scalar units implies that the permutations $\mathbf e_\mu\mapsto\mathbf e_{u\mu}$ commute with $\rho_L$.
They therefore act pointwise on $M_{k,\rho_L}$ and $S_{k,\rho_L}$, and we use a superscript $\mathcal U$ to denote the invariant subspaces.

For $\gamma=\left(\begin{smallmatrix}a&b\\c&d\end{smallmatrix}\right)\in\mathrm U(1,1)(\mathbb Z)$, use the fractional linear action $\gamma\tau=(a\tau+b)/(c\tau+d)$ and the automorphy factor $j(\gamma,\tau)=c\tau+d$.
Define the slash operator and the spaces $M_{k,\rho_L^{\mathrm U}}$ and $S_{k,\rho_L^{\mathrm U}}$ as in Definition~\ref{def:vectorValuedModularForms}, with $\mathrm U(1,1)(\mathbb Z)$ in place of $\operatorname{SL}_2(\mathbb Z)$.
The two groups have the same effective action on $\mathbb H$, so the cusp condition is again imposed at $\infty$.
For $f=\sum_{\mu\in A_L}f_\mu\mathbf e_\mu\in M_{n+1,\rho_L}$, the transformation law under a scalar unit $uI_2$ is equivalent to $f_{u\mu}=f_\mu$, since $uI_2$ fixes $\mathbb H$ pointwise and $j(uI_2,\tau)=u$.
The group decomposition therefore gives $M_{n+1,\rho_L^{\mathrm U}}=M_{n+1,\rho_L}^{\mathcal U}$ and
\begin{equation}
S_{n+1,\rho_L^{\mathrm U}}=S_{n+1,\rho_L}^{\mathcal U}.
\label{eq:unitaryCuspSpace}
\end{equation}
In particular, the Fourier coefficients of $g\in S_{n+1,\rho_L^{\mathrm U}}$ satisfy $b_g(m,u\mu)=b_g(m,\mu)$ for every $u\in\mathcal U$.
\subsection{Harmonic weak Maass forms}
\label{subsec:harmonicWeakMaassForms}

For $\tau=u+iv\in\mathbb H$ and $k\in\mathbb Z$, the weight $k$ hyperbolic Laplacian is
$$
\Delta_k=-v^2\left(\frac{\partial^2}{\partial u^2} +\frac{\partial^2}{\partial v^2}\right)+ikv\left( \frac{\partial}{\partial u}+i\frac{\partial}{\partial v}\right).
$$

\begin{definition}
\label{def:harmonicWeakMaassForms}
For $(W,\rho)$ as above, a harmonic weak Maass form of weight $k$ and type $\rho$ is a smooth function $f:\mathbb H\to W$ satisfying
\begin{itemize}
    \item $f\vert_{k,\rho}\gamma=f$ for all $\gamma\in\operatorname{SL}_2(\mathbb Z)$,
    \item $\Delta_kf=0$,
    \item $f(u+iv)=O(e^{Cv})$ as $v\to\infty$ for some $C>0$, uniformly in $u$.
\end{itemize}
We denote the space of such forms by $H_{k,\rho}$.
\end{definition}

The invariant Hermitian form identifies $\overline W$ with $W^\vee$.
As in \cite[Proposition 3.2]{BFtwolift}, the operator
$$
\xi_{k,\rho}(f)(\tau)=2iv^k\overline{\frac{\partial f}{\partial\overline\tau}}
$$
defines a conjugate linear map $H_{k,\rho}\to M^!_{2-k,\rho^\vee}$.
We write $\xi_{1-n}:H_{1-n,\overline{\rho_L}}\longrightarrow M^!_{n+1,\rho_L}$ for its specialization to weight $1-n$ and type $\overline{\rho_L}$, and set $H^+_{1-n,\overline{\rho_L}}=\xi_{1-n}^{-1}(S_{n+1,\rho_L})$.

Every $f\in H^+_{1-n,\overline{\rho_L}}$ has a unique decomposition $f=f^++f^-$ into its holomorphic and nonholomorphic parts, with Fourier expansions
$$
\begin{aligned}
f^+(\tau)
&=\sum_{\mu\in A_L}\sum_{\substack{m\gg-\infty\\m+\mathbb Z=-q_L(\mu)}}c_f^+(m,\mu)e(m\tau)\mathbf e_\mu,\\
f^-(\tau)
&=\sum_{\mu\in A_L}\sum_{\substack{m>0\\
-m+\mathbb Z=-q_L(\mu)}}c_f^-(-m,\mu)\Gamma(n,4\pi mv)e(-m\tau)\mathbf e_\mu.
\end{aligned}
$$
Here $\Gamma(a,x)=\int_x^\infty t^{a-1}e^{-t}\,dt$ denotes the incomplete gamma function.
Because $\xi_{1-n}(f)$ is cuspidal, $f^-$ has no $m=0$ term and hence decays exponentially as $v\to\infty$, uniformly in $u$.
If $f\in M^!_{1-n,\overline{\rho_L}}$, then $f^-=0$ and $c_f^+=c_f$.

The principal part of $f\in H^+_{1-n,\overline{\rho_L}}$ is the finite Fourier polynomial
$$
P_f(\tau)=\sum_{\mu\in A_L}\sum_{\substack{m<0\\m+\mathbb Z=-q_L(\mu)}}c_f^+(m,\mu)e(m\tau)\mathbf e_\mu.
$$
Unlike \cite[(3.6)]{BFtwolift}, we omit the constant term. We say that $P_f$ is integral if $c_f^+(m,\mu)\in\mathbb Z$ for all $m<0$ and $\mu\in A_L$.

Since $\rho_L(-I)\mathbf e_\mu=\overline{\rho_L}(-I)\mathbf e_\mu=(-1)^{n-1}\mathbf e_{-\mu}$, invariance under $-I$ forces the Fourier coefficients of $f\in H^+_{1-n,\overline{\rho_L}}$ and $g\in S_{n+1,\rho_L}$ to satisfy the symmetry relations $c_f^+(m,\mu)=c_f^+(m,-\mu)$ and $b_g(m,\mu)=b_g(m,-\mu)$.
In particular, every principal part arising from $H^+_{1-n,\overline{\rho_L}}$ is symmetric, in the sense that its coefficients satisfy $c(m,\mu)=c(m,-\mu)$.

\begin{proposition}[{\cite[Proposition 3.2]{BFtwolift}}]
\label{prop:xiKernel}
For $f\in H^+_{1-n,\overline{\rho_L}}$, one has $\xi_{1-n}(f)=0$ if and only if $f\in M^!_{1-n,\overline{\rho_L}}$.
\end{proposition}

\subsection{Unitary Siegel--Weil formula and the Eisenstein series}
\label{subsec:siegelWeilEisenstein}

Let $(V',h')$ be a Hermitian space over $K$ of dimension $d\geq3$ and signature $(d-1,1)$.
We retain the notation and normalizations of Subsection \ref{subsec:weilRepresentations} and identify $\operatorname{SL}_2$ with its image in $G=\operatorname{U}(1,1)$ under $\jmath$.
Let $P_G=N_GM_G$ be the standard upper triangular parabolic subgroup of $G$, where $N_G$ is the unipotent radical and $M_G$ is the diagonal Levi subgroup.
Put $K_G=K_{G,\infty}\prod_pK_{G,p}$, where $K_{G,p}=G(\mathbb Q_p)\cap\operatorname{GL}_2(\mathcal O_{K,p})$ and $K_{G,\infty}=G(\mathbb R)\cap\operatorname{U}(2)$.

For $\varphi\in\mathcal S(V'(\mathbb A_K))$, $g\in G(\mathbb A)$, and $h\in H'(\mathbb A)$, define
$$
\theta(g,h;\varphi)=\sum_{x\in V'(K)}(\omega_{V'}(g,h)\varphi)(x).
$$
Equip $[H']:=H'(\mathbb Q)\backslash H'(\mathbb A)$ with the invariant measure $dh$ normalized by $\operatorname{vol}([H'],dh)=1$.
We set
\begin{equation}
I(g;\varphi):=\int_{[H']}\theta(g,h;\varphi)\,dh.
\label{eq:adelicThetaIntegral}
\end{equation}
We add the subscript $V'$ to $\theta$ and $I$ whenever the underlying Hermitian space needs to be specified.

For $s\in\mathbb C$, the normalized principal series $I(s,\chi')$ consists of smooth $K_G$-finite functions $\Phi(\,\cdot\,,s)$ on $G(\mathbb A)$ satisfying
$$
\Phi(n(b)m(a)g,s)=\chi'(a)|a|_{\mathbb A_K}^{s+1/2}\Phi(g,s).
$$
A family $\Phi=(\Phi(\,\cdot\,,s))_{s\in\mathbb C}$ with $\Phi(\,\cdot\,,s)\in I(s,\chi')$ is called a standard section if $\Phi(k,s)$ is independent of $s$ for every $k\in K_G$.
For such a section $\Phi$, the associated Eisenstein series is
$$
E(g,s,\Phi)=\sum_{\gamma\in P_G(\mathbb Q)\backslash G(\mathbb Q)}\Phi(\gamma g,s),
\qquad \operatorname{Re}(s)>\frac12.
$$
Setting $s_0=(d-1)/2$, we associate to each $K_G$-finite Schwartz function $\varphi$ the standard Siegel--Weil section
$$
\lambda(\varphi)(g,s)=|a(g)|_{\mathbb A_K}^{s-s_0}(\omega_{V'}(g,1)\varphi)(0),
$$
where $g=n(b)m(a(g))k$ is an Iwasawa decomposition.

The convergent unitary Siegel--Weil formula takes the following form.
\begin{theorem}[{\cite[Theorem 2]{YamanaunitarySW}}]
\label{thm:convergentUnitarySiegelWeil}
Let $r_{V'}$ denote the Witt index of $V'$ over $K$, and assume that $d-r_{V'}>1$.
For every $K_G$-finite Schwartz function $\varphi\in\mathcal S(V'(\mathbb A_K))$ and every $g\in G(\mathbb A)$, the theta integral $I(g;\varphi)$ converges absolutely.
Moreover, the Eisenstein series $E(g,s,\lambda(\varphi))$ is holomorphic at $s=s_0$ and satisfies
$$
E(g,s_0,\lambda(\varphi))=c_{V'}I(g;\varphi),
$$
where $c_{V'}$ is a nonzero constant independent of $g$ and $\varphi$.
\end{theorem}

Let $J\subset V'$ be a trace even Hermitian lattice, and set $\Gamma_\infty=P_G(\mathbb Q)\cap\operatorname{SL}_2(\mathbb Z) =\{\pm T^m:m\in\mathbb Z\}$. 
For $\tau=u+iv\in\mathbb H$ and $\operatorname{Re}(s)>1/2$, the $\mathbb C[A_J]$-valued Eisenstein series is defined by
\begin{equation}
E_J(\tau,s,d):=\sum_{\gamma\in\Gamma_\infty\backslash\operatorname{SL}_2(\mathbb Z)}\operatorname{Im}(\gamma\tau)^{s+(1-d)/2}(c\tau+d_\gamma)^{-d}\,\rho_J(\gamma)^{-1}\mathbf e_0,
\label{eq:vectorValuedEisenstein}
\end{equation}
where $(c,d_\gamma)$ is the lower row of $\gamma$.
Since $\sigma_J=2d-4$, we have $\rho_J(T)\mathbf e_0=\mathbf e_0$ and $\rho_J(-I)\mathbf e_0=(-1)^d\mathbf e_0$, which ensures that each summand in \eqref{eq:vectorValuedEisenstein} is independent of the choice of coset representative.
This series converges absolutely and locally uniformly on $\mathbb H\times\{s\in\mathbb C:\operatorname{Re}(s)>1/2\}$, defining a holomorphic function of $s$ in this half-plane.

Let $\Phi_\infty$ be a normalized archimedean standard section of weight $d$. 
Explicitly, $\Phi_\infty(1,s)=1$ and $\Phi_\infty(g\jmath(k),s)=\Phi_\infty(g,s)(a+ib)^d$ for $k=\left(\begin{smallmatrix}a&b\\-b&a\end{smallmatrix}\right)$.
For $\mu\in A_J$, let $\Phi_{\mu,f}$ be the finite standard section
attached to $\phi_\mu^J$ and put
$\Phi_\mu=\Phi_\infty\otimes\Phi_{\mu,f}$.

\begin{proposition}
\label{prop:adelicClassicalEisenstein}
For $\tau=u+iv\in\mathbb H$ and $\operatorname{Re}(s)>1/2$, one has
$$
v^{-d/2}\sum_{\mu\in A_J} E(\jmath(g_\tau),s,\Phi_\mu)\mathbf e_\mu=E_J(\tau,s,d).
$$
\end{proposition}

\begin{proof}
For $\gamma\in\operatorname{SL}_2(\mathbb Z)$, the definition of the finite standard sections and Proposition \ref{prop:finiteSchwartzModel} imply
$$
\sum_{\mu\in A_J}\Phi_{\mu,f}(\jmath(\gamma_f),s)\mathbf e_\mu=\rho_J(\gamma)^{-1}\mathbf e_0.
$$
Applying \cite[Lemma 2.2.1]{BruinierHoward} to each $\Phi_\mu$, with $\ell=d$ and $j(g_\tau,i)=v^{-1/2}$ in their notation, then yields the desired identity.
\end{proof}

\section{Borcherds products and the weak converse theorem}
\label{sec:borcherdsConverse}
Throughout this section, we retain the Hermitian lattice $L\subset V$ of signature $(n,1)$ from Section \ref{sec:unitary shimura varieties}.
Here we formulate the weak converse theorem and explain how the injectivity of the unitary Kudla--Millson lift implies the converse theorem.

\subsection{Orthogonal and unitary Borcherds products}
\label{subsec:borcherdsProducts}

A Borcherds product is a meromorphic automorphic form constructed from a weakly holomorphic vector-valued modular form.
Its divisor is a linear combination of Heegner divisors determined by the negative index Fourier coefficients of the vector-valued modular form.

Let $M$ be an even lattice of signature $(r,2)$ with $r\geq3$.
When $r$ is odd, the modular form spaces in the following theorem are understood in the usual half-integral weight setting for the metaplectic double cover $\operatorname{Mp}_2(\mathbb Z)$ and the Weil representation attached to $(A_M,-q_M)$.
Borcherds's automorphic product theorem states the following.

\begin{theorem}[{\cite[Theorem 13.3]{borcherds1996automorphic}; see also \cite[Theorem 1.1]{Bruinier14Converse}}]
\label{thm:orthogonalBorcherdsProduct}
Let $f\in M^{!}_{1-r/2,\overline{\rho_M}}$ be a weakly holomorphic modular form whose Fourier coefficients $c_f(m,{\mu})$ are integral for $m<0$.
Then there exists a meromorphic automorphic form
$\widetilde\Psi_M(\,\cdot\,,f)$ for $\Gamma_M^{\mathrm O}$ of rational weight $c_f(0,0)/2$ with a multiplier system of
finite order. 
Its divisor on $\mathcal D_M$ is given by
$$
\operatorname{div}_{\mathcal D_M}\bigl(\widetilde\Psi_M(\,\cdot\,,f)\bigr)=\frac12 \sum_{{\mu}\in A_M}\sum_{\substack{m<0\\m+\mathbb Z=-q_M({\mu})}}c_f(m,{\mu})\widetilde H_{\mathrm{O}}(-m,{\mu}).
$$
The factor $1/2$ accounts for the equal contributions of opposite vectors in the Heegner divisor sums.
\end{theorem}
The meromorphic automorphic form $\widetilde\Psi_M(\,\cdot\,,f)$ admits an infinite product expansion in suitable coordinates near each zero-dimensional cusp and is therefore called the \emph{Borcherds product} associated with $f$.

For the remainder of this section, we assume that $n\geq2$.
By restricting the orthogonal Borcherds product along the trace embedding $\iota$, Hofmann obtained the following unitary analogue.
\begin{theorem}[{\cite[Theorem 8.1(1)--(2)]{hofmann2014borcherds}}]
\label{thm:unitaryBorcherdsProduct}
Let $f\in M^!_{1-n,\overline{\rho_L}}$ be a weakly holomorphic modular form whose Fourier coefficients $c_f(m,{\mu})$ are integral for $m<0$.
Then there exists a meromorphic automorphic form $\widetilde\Psi_{\mathrm{U}}(\,\cdot\,,f)$ for $\Gamma_L^{\mathrm U}$ of rational weight $c_f(0,0)/2$ with a multiplier system of finite order.
Its divisor on $\mathcal D_{V}$ is
$$
\operatorname{div}_{\mathcal D_{V}}\bigl(\widetilde\Psi_{\mathrm{U}}(\,\cdot\,,f)\bigr)=\frac12\sum_{{\mu}\in A_L}\sum_{\substack{m<0\\m+\mathbb Z=-q_L({\mu})}} c_f(m,{\mu})\widetilde H_{\mathrm{U}}(-m,{\mu}).
$$
\end{theorem}

\subsection{The unitary singular theta lift}
\label{subsec:singularThetaLift}
The logarithmic Petersson norm of a Borcherds product can be expressed in terms of a regularized theta integral, called the singular theta lift.
We construct it from the unitary Siegel theta function following \cite{FunkeHofmannGreen} and compare the resulting lift with Bruinier's orthogonal lift.

For $z\in\mathcal D_{V}$ and $x\in V_{\mathbb C}$, write $x=x_z+x_{z^\perp}$ for the $h$-orthogonal decomposition.
Define the Gaussian
$
\varphi_0(x,z)=\exp\!\left(-2\pi\bigl(h(x_{z^\perp},x_{z^\perp})-h(x_z,x_z)\bigr)\right).
$
This normalization matches the Weil representation fixed in Subsection \ref{subsec:weilRepresentations} and coincides with the Gaussian of \cite[Section 3.1]{FunkeHofmannGreen} evaluated at $\sqrt2\,x$.
For $\tau=u+iv\in\mathbb H$, put
$$
\varphi_0(x,\tau,z):=v^{-(n-1)/2}\bigl(\omega_{V,\infty}(\jmath(g_\tau),1)\varphi_0(\,\cdot\,,z)\bigr)(x)=v\,e\bigl(uQ(x)\bigr)\varphi_0(\sqrt v\,x,z).
$$

\begin{definition}
\label{def:unitarySiegelTheta}
For $(\tau,z)\in\mathbb H\times\mathcal D_{V}$ and ${\mu}\in A_L$, define the unitary scalar Siegel theta functions and the vector-valued Siegel theta function by
$$
\begin{aligned}
\theta_{L,{\mu}}^{\mathrm{U}}(\tau,z)
&:=\sum_{\substack{\lambda\in L^\vee\\\lambda+L={\mu}}}\overline{\varphi_0(\lambda,\tau,z)},\\
\Theta_L^{\mathrm{U},\vee}(\tau,z)&:=\sum_{{\mu}\in A_L}\overline{\theta_{L,{\mu}}^{\mathrm{U}}(\tau,z)}\,\mathbf e_{\mu}^\vee.
\end{aligned}
$$
\end{definition}
Since $Q(x)=h(x,x)$, this definition yields
\begin{equation}
\theta_{L,{\mu}}^{\mathrm{U}}(\tau,z)=v\sum_{\substack{\lambda\in L^\vee\\\lambda+L={\mu}}}e\!\left(-\overline\tau Q(\lambda)-(\tau-\overline\tau)h(\lambda_z,\lambda_z)\right).
\label{eq:scalarSiegelTheta}
\end{equation}
The theta transformation law, together with the weight $n-1$ of the Gaussian and Proposition \ref{prop:finiteSchwartzModel}, shows that $\Theta_L^{\mathrm{U},\vee}(\tau,z)$ transforms with weight $n-1$ and type $(\overline{\rho_L})^\vee$ in $\tau$.
Thus, for $f\in H^+_{1-n,\overline{\rho_L}}$, the function $\{f(\tau),\Theta_L^{\mathrm{U},\vee}(\tau,z)\}$ is $\operatorname{SL}_2(\mathbb Z)$-invariant.

To compare with the orthogonal kernel, take $M=L_{\mathbb Z}$ with quadratic form $Q$.
For $Z\in\mathcal D_M$ and $\lambda\in M\otimes_{\mathbb Z}\mathbb R$, write $\lambda=\lambda_Z+\lambda_{Z^\perp}$ for the orthogonal decomposition, and set
$$
\begin{aligned}
\theta_{M,{\mu}}^{\mathrm{O}}(\tau,Z)
&=v\sum_{\substack{\lambda\in M^\vee\\\lambda+M={\mu}}}e\!\left(-\tau Q(\lambda_Z)-\overline\tau Q(\lambda_{Z^\perp})\right),\\
\Theta_M^{\mathrm{O},\vee}(\tau,Z)
&=\sum_{{\mu}\in A_L} \overline{\theta_{M,{\mu}}^{\mathrm{O}}(\tau,Z)}\,\mathbf e_{\mu}^\vee.
\end{aligned}
$$
For $Z=\iota(z)$, the orthogonal projections onto $Z$ and $Z^\perp$ coincide with the Hermitian projections onto $z$ and $z^\perp$, so that $Q(\lambda_Z)=h(\lambda_z,\lambda_z)$.
Comparing the preceding formula with \eqref{eq:scalarSiegelTheta} then yields
$$
\theta_{L,{\mu}}^{\mathrm{U}}(\tau,z)=\theta_{M,{\mu}}^{\mathrm{O}}(\tau,\iota(z)),
\qquad
\Theta_L^{\mathrm{U},\vee}(\tau,z)=\Theta_M^{\mathrm{O},\vee}(\tau,\iota(z)).
$$

To integrate this theta kernel against $f$, let $\mathcal F$ be the standard fundamental domain for $\operatorname{SL}_2(\mathbb Z)$, truncated at height $T$ by $\mathcal F_T=\{\tau=u+iv\in\mathcal F:v\leq T\}$, and let $d\mu(\tau)=du\,dv/v^2$ be the invariant hyperbolic measure on $\mathbb H$.
Writing $\operatorname{CT}_{w=0}$ for the constant term in the Laurent expansion at $w=0$, we define the regularized integral of a smooth function $A:\mathbb H\to\mathbb C$ by
$$\int_{\mathcal F}^{\mathrm{reg}}A(\tau)\,d\mu(\tau):=\operatorname{CT}_{w=0}\left[\lim_{T\to\infty}\int_{\mathcal F_T}A(\tau)v^{-w}\,d\mu(\tau)\right]$$
whenever the limit exists for $\operatorname{Re}(w)\gg 0$ and admits a meromorphic continuation to $w=0$.
For $f\in H^+_{1-n,\overline{\rho_L}}$, define the $\mathbb C$-divisor associated with $f$ on $\mathcal D_M$ by
$$
\widetilde H_{\mathrm{O}}(f)=\frac12\sum_{{\mu}\in A_L}\sum_{\substack{m<0\\m+\mathbb Z=-q_L({\mu})}} c_f^+(m,{\mu})\widetilde H_{\mathrm{O}}(-m,{\mu}).
$$
Away from the hyperplanes arising from nonzero principal part coefficients, define the orthogonal singular theta lift of $f$ by
$$
\widetilde\Phi_M^{\mathrm{O}}(Z,f)=\int_{\mathcal F}^{\mathrm{reg}}\{f(\tau),\Theta_M^{\mathrm{O},\vee}(\tau,Z)\}\,d\mu(\tau).
$$
The lift is $\Gamma_M^{\mathrm O}$-invariant.
By \cite[Proposition 5.4]{BFtwolift}, it extends real analytically to $\mathcal D_M\setminus|\widetilde H_{\mathrm O}(f)|$ and has logarithmic singularities along $\widetilde H_{\mathrm O}(f)$.

On $\mathcal D_{V}$, define the divisor associated with $f$ by
$$
\widetilde H_{\mathrm{U}}(f)=\frac12\sum_{{\mu}\in A_L}\sum_{\substack{m<0\\m+\mathbb Z=-q_L({\mu})}}c_f^+(m,{\mu})\widetilde H_{\mathrm{U}}(-m,{\mu}),
$$
and let $H_{\mathrm{U}}(f)$ denote its descent to $X_L$.
Away from the hyperplanes arising from nonzero principal part coefficients, define the unitary singular theta lift by
$$
\widetilde\Phi_L^{\mathrm{U}}(z,f)=\int_{\mathcal F}^{\mathrm{reg}}\{f(\tau),\Theta_L^{\mathrm{U},\vee}(\tau,z)\}\,d\mu(\tau).
$$
It is $\Gamma_L^{\mathrm U}$-invariant and descends to $\Phi_L^{\mathrm{U}}(\,\cdot\,,f)$ on $X_L$.
For a finite index subgroup $\Gamma\subseteq\Gamma_L^{\mathrm U}$, we use the same notation for the descents of $\widetilde H_{\mathrm U}(f)$ and $\widetilde\Phi_L^{\mathrm U}(\,\cdot\,,f)$ to $X_\Gamma$.
\begin{proposition}
\label{prop:singularThetaRegularity}
Let $f\in H^+_{1-n,\overline{\rho_L}}$.
The lift $\widetilde\Phi_L^{\mathrm U}(\,\cdot\,,f)$ extends smoothly to $\mathcal D_V\setminus|\widetilde H_{\mathrm U}(f)|$ and defines a locally integrable function on $\mathcal D_V$.
Moreover, let $W\subset\mathcal D_{V}$ be an open neighborhood and let $F_W$ be a nonzero meromorphic function on $W$ satisfying
$\operatorname{div}_W(F_W)=\widetilde H_{\mathrm{U}}(f)|_W$.
Then $\widetilde\Phi_L^{\mathrm{U}}(\,\cdot\,,f)+2\log|F_W|^2$ extends smoothly from $W\setminus|\widetilde H_{\mathrm{U}}(f)|$ to $W$.
\end{proposition}

\begin{proof}
Fix $z_0\in\mathcal D_V$ and choose a neighborhood $W$ small enough that every contributing hyperplane meeting $W$ passes through $z_0$.
In our notation, \cite[Theorem 4.1]{BruinierHowardYang2015} yields
$$
\widetilde\Phi_L^{\mathrm{U}}(z,f)=H_W(z)+\sum_{{\mu}\in A_L} \sum_{\substack{m<0\\m+\mathbb Z=-q_L({\mu})}}c_f^+(m,{\mu})\sum_{\substack{\lambda\in({\mu}+L)\cap z_0^\perp\\Q(\lambda)=-m}}\Gamma\bigl(0,-4\pi h(\lambda_z,\lambda_z)\bigr)
$$
for $z\in W$ away from the contributing hyperplanes, where $H_W$ is smooth on $W$.

Choose a nowhere-vanishing holomorphic section $s$ of the tautological line bundle over $W$, so that $-h(\lambda_z,\lambda_z)= |h(s(z),\lambda)|^2/\bigl(-h(s(z),s(z))\bigr)$.
Since $\Gamma(0,x)+\log x$ extends smoothly to $x=0$, each summand differs from $-c_f^+(-Q(\lambda),\lambda+L)\log|h(s(z),\lambda)|^2$ by a smooth function on $W$.
Because $h(s(z),\lambda)$ vanishes to order one along the hyperplane $\mathcal D_\lambda:=\{z\in\mathcal D_{V}:z\subset\lambda^\perp\}=\mathcal D_{V}\cap\mathbb P(\lambda^\perp)$, the definition of $\widetilde H_{\mathrm{U}}(f)$ implies that
$$
\sum_\lambda c_f^+(-Q(\lambda),\lambda+L)\,\operatorname{div}_W\bigl(h(s(\,\cdot\,),\lambda)\bigr)=2\widetilde H_{\mathrm{U}}(f)|_W,
$$
where the sum is over the same vectors as above.

If $\operatorname{div}_W(F_W)=\widetilde H_{\mathrm{U}}(f)|_W$, the logarithmic terms therefore combine to $-2\log|F_W|^2$ up to a smooth function.
Hence $\widetilde\Phi_L^{\mathrm{U}}(\,\cdot\,,f)+2\log|F_W|^2$ extends smoothly to $W$.

The expansion and divisor identity above show that the regularized theta lift $\widetilde\Phi_L^{\mathrm{U}}(\,\cdot\,,f)$ is locally integrable on $W$ and has a smooth extension on $W\setminus|\widetilde H_{\mathrm{U}}(f)|$.
Since these smooth extensions coincide with the regularized theta integral away from the locally finite collection of hyperplanes $\mathcal D_\lambda$, they agree on overlaps by continuity and thus glue to a smooth extension on $\mathcal D_{V}\setminus|\widetilde H_{\mathrm{U}}(f)|$.
\end{proof}

Let $\mathcal L$ be the tautological line bundle over $\mathcal D_{V}$, whose fiber at $z\in\mathcal D_{V}$ is the line $z$ equipped with the Hermitian metric $\lVert v\rVert_z^2=-h(v,v)$ for $v\in z$.
We regard a meromorphic automorphic form $F$ of integral weight $r$ with a finite order multiplier system as a meromorphic section of $\mathcal L^{\otimes r}$ on $\mathcal D_V$.
Locally, write $F=F_s s^{\otimes r}$, where $s$ is a nowhere-vanishing holomorphic section of $\mathcal L$.
The metric on $\mathcal L$ induces the Petersson norm $\lVert F\rVert=|F_s|\lVert s\rVert^r$, which does not depend on the choice of $s$. 
When the weight $r$ is rational, choosing a positive integer $N$ with $Nr\in\mathbb Z$ allows us to define $\lVert F\rVert=\lVert F^N\rVert^{1/N}$ independently of $N$.
The Petersson norm is $\Gamma$-invariant since the metric is $\Gamma$-invariant and the multiplier has absolute value one.
We also write $\lVert F\rVert$ for its descent to $X_\Gamma$.

For $f\in M^!_{1-n,\overline{\rho_L}}$ with integral principal part, Hofmann's formulas \cite[(27)--(28)]{hofmann2014borcherds}, expressed in our normalization, yield a real constant $C_f^{\mathrm{U}}$ such that
$$
\Phi_L^{\mathrm{U}}(\,\cdot\,,f)=-4\log\lVert\widetilde\Psi_{\mathrm{U}}(\,\cdot\,,f)\rVert+C_f^{\mathrm{U}}
$$
on $X_\Gamma\setminus|H_{\mathrm{U}}(f)|$.
\subsection{The Kudla--Millson lift}
\label{subsec:unitaryKMLift}
We use $d^c=(\partial-\overline\partial)/(4\pi i)$, so that $dd^c=(i/2\pi)\partial\overline\partial$.
Let $(V',h')$ be a Hermitian space over $K$ of dimension $d\geq3$ and signature $(d-1,1)$, and put $H'=\mathrm U(V')$.
Let $\mathcal D'$ be the domain of negative lines in $V'\otimes_{K,\sigma}\mathbb C$.
The metric $\lVert v\rVert_z^2=-h'(v,v)$ on the tautological line bundle over $\mathcal D'$ defines an $H'(\mathbb R)$-invariant K\"ahler form $\Omega'$, given locally by $\Omega'=-dd^c\log\lVert s\rVert^2$ for any nowhere-vanishing holomorphic section $s$.

We use the Kudla--Millson Schwartz form of \cite[Section 3.2]{FunkeHofmannGreen} for $(p,q)=(d-1,1)$ with $(x,y)_{\mathrm{FH}}=h'(y,x)$, and write $\varphi_{\mathrm{KM},V'}(x,z)$ for this form evaluated at $\sqrt2\,x$ to match the Weil representation fixed in Subsection \ref{subsec:weilRepresentations}.
The form $\varphi_{\mathrm{KM},V'}(x,z)$ is Schwartz in $x$ and takes values in smooth closed $(1,1)$-forms in $z\in\mathcal D'$.
When the underlying Hermitian space is clear, we omit the subscript $V'$ and view $\varphi_{\mathrm{KM}}(x)=\varphi_{\mathrm{KM}}(x,\cdot)$ as a differential form on $\mathcal D'$ satisfying $\varphi_{\mathrm{KM}}(0)=i\,\Omega'$.
Define the scalarization $\widetilde\varphi_{\mathrm{KM}}(x,z)$ by
\begin{equation}
\varphi_{\mathrm{KM}}(x,z)\wedge\Omega'(z)^{d-2}
=\widetilde\varphi_{\mathrm{KM}}(x,z)\,\Omega'(z)^{d-1}.
\label{eq:kmScalarization}
\end{equation}
By the $H'(\mathbb R)$-equivariance of the Kudla--Millson form, we have $\widetilde\varphi_{\mathrm{KM}}(h_\infty x,h_\infty z)=\widetilde\varphi_{\mathrm{KM}}(x,z)$
for $h_\infty\in H'(\mathbb R)$.
Fixing a base point in $\mathcal D'$ for each such Hermitian space, we write $\widetilde\varphi_{\mathrm{KM},V'}$ without the variables for the scalar Schwartz function on $V'(\mathbb R)$ obtained by evaluating at this base point.

Let $\varphi_f\in\mathcal S(V'(\mathbb A_{K,f}))$ and let $K_f\subset H'(\mathbb A_f)$ be a compact open subgroup fixing $\varphi_f$.
On the adelic quotient $X_{K_f}$ of Section \ref{sec:unitary shimura varieties} associated with $(V',h')$, each $g=(g_\infty,g_f)\in G(\mathbb A)$ defines a $(1,1)$-form $\theta(g;\varphi_{\mathrm{KM}}\otimes\varphi_f)$ whose pullback at $(z,h_f)\in\mathcal D'\times H'(\mathbb A_f)$ is given by
$$
\sum_{x\in V'(K)}\bigl[\omega_{V',\infty}(g_\infty,1)\varphi_{\mathrm{KM}}\bigr](x_\infty,z) \bigl[\omega_{V',f}(g_f,h_f)\varphi_f\bigr](x_f).
$$
After pullback to $H'(\mathbb A)$ using the fixed base point, one has
$$
\theta(g;\varphi_{\mathrm{KM}}\otimes\varphi_f)\wedge\Omega'^{d-2}=\theta(g,h;\widetilde\varphi_{\mathrm{KM}}\otimes\varphi_f)\Omega'^{d-1}.
$$

The scalarized function has weight $d$ under $\operatorname{SO}(2)$ by \cite[Theorem 3.1(1)]{FunkeHofmannGreen} and satisfies $\widetilde\varphi_{\mathrm{KM}}(0)=i$.
Hence, for the sections associated with a lattice $J\subset V'$ in Subsection \ref{subsec:siegelWeilEisenstein}, we have $\lambda(\widetilde\varphi_{\mathrm{KM}}\otimes\phi_\mu^J)=i\,\Phi_\mu$.

Specializing to $V'=V$, we have $d=n+1$ and $\mathcal D'=\mathcal D_{V}$, and we write $\Omega=\Omega'$.
For $\tau=u+iv$, put $\varphi_{\mathrm{KM}}(x,\tau,z)=e(uQ(x))\varphi_{\mathrm{KM}}(\sqrt v\,x,z)$ and define
$$
\Theta_L^{\mathrm{KM}}(\tau,z)
=\sum_{\mu\in A_L}\sum_{\substack{\lambda\in L^\vee\\\lambda+L=\mu}}
\varphi_{\mathrm{KM}}(\lambda,\tau,z)\mathbf e_\mu.
$$
On the canonical component $X_L\subset X_{K_L}$, this series agrees with
$$
v^{-(n+1)/2}\sum_{\mu\in A_L}
\left.\theta\bigl(\jmath(g_\tau);\varphi_{\mathrm{KM}}\otimes\phi_\mu^L\bigr)\right|_{X_L}\mathbf e_\mu.
$$
Under the coefficient module identification of Subsection \ref{subsec:weilRepresentations}, this agrees with the Kudla--Millson theta series in \cite[Section 5.2]{FunkeHofmannGreen}.
Proposition \ref{prop:finiteSchwartzModel} and \cite[Theorem 3.1(1)--(2)]{FunkeHofmannGreen} show that $\Theta_L^{\mathrm{KM}}$ has weight $n+1$ and type $\rho_L$ in $\tau$, with values in smooth closed $(1,1)$-forms on $\mathcal D_{V}$.

Let $\mathcal Z^{1,1}(X_{\Gamma})$ denote the space of smooth closed $(1,1)$-forms on $X_{\Gamma}$.
For $g\in S_{n+1,\rho_L}$, define
$$
\Lambda_{\mathrm{U}}(g)(z):=\int_{\mathcal F} \left\langle\Theta_L^{\mathrm{KM}}(\tau,z),g(\tau)\right\rangle_{A_L} v^{n+1}\,\frac{du\,dv}{v^2}.
$$
The equivariance of the Kudla--Millson Schwartz form and the trivial action of $\Gamma_L^{\mathrm U}$ on $A_L$ show that $\Theta_L^{\mathrm{KM}}(\tau,\cdot)$ is $\Gamma_L^{\mathrm U}$-invariant and hence descends to $X_\Gamma$.
Since the components of the theta function are smooth closed forms and the integral converges absolutely, $\Lambda_{\mathrm{U}}$ defines a conjugate linear map into $\mathcal Z^{1,1}(X_{\Gamma})$.

\begin{remark}
For $\varepsilon\in\mathcal U$, the scalar $\varepsilon \operatorname{id}_V$ fixes $\mathcal D_V$ pointwise, so the equivariance of the Kudla--Millson Schwartz form gives
$\varphi_{\mathrm{KM}}(\varepsilon x,\tau,z)=\varphi_{\mathrm{KM}}(x,\tau,z)$.
Reindexing the lattice sum via $x\mapsto\varepsilon x$ shows that $\Theta_L^{\mathrm{KM}}$ is $\mathcal U$-invariant.
Together with its $\operatorname{SL}_2(\mathbb Z)$ transformation law, this implies that $\Theta_L^{\mathrm{KM}}$ satisfies the transformation law of weight $n+1$ and type $\rho_L^{\mathrm U}$ for $\mathrm U(1,1)(\mathbb Z)$.

Since $\mathrm U(1,1)(\mathbb Z)$ and $\operatorname{SL}_2(\mathbb Z)$ have the same effective action on $\mathbb H$, we use the same fundamental domain $\mathcal F$ and measure for both groups.
By \eqref{eq:unitaryCuspSpace}, restricting $\Lambda_{\mathrm U}$ to $S_{n+1,\rho_L^{\mathrm U}}=S_{n+1,\rho_L}^{\mathcal U}$ gives the Kudla--Millson lift for $\mathrm U(1,1)(\mathbb Z)$, denoted by the same symbol:
$$
\Lambda_{\mathrm U}:S_{n+1,\rho_L^{\mathrm U}}
\longrightarrow\mathcal Z^{1,1}(X_\Gamma).
$$
\end{remark}

The singular theta lift satisfies the following curvature identity.

\begin{proposition}[{cf. \cite[Theorem 5.6(1)]{FunkeHofmannGreen}}]
\label{prop:thetaLiftCurvature}
Let $f\in H^+_{1-n,\overline{\rho_L}}$, and put $g=\xi_{1-n}(f)$.
The form $dd^c\Phi_L^{\mathrm{U}}(\,\cdot\,,f)$ on $X_{\Gamma}\setminus|H_{\mathrm{U}}(f)|$ extends uniquely to a smooth closed $(1,1)$-form on $X_{\Gamma}$.
Denoting this extension by the same symbol, one has
$$
dd^c\Phi_L^{\mathrm{U}}(\,\cdot\,,f) =i\Lambda_{\mathrm{U}}(g)+c_f^+(0,0)\Omega.
$$
\end{proposition}

\begin{proof}
Although Funke and Hofmann assume throughout their paper that the Hermitian form is $\mathcal O_K$-valued on the lattice, their proof of \cite[Theorem 5.6]{FunkeHofmannGreen} remains valid for the trace even Hermitian lattices considered here.
For the reader's convenience, we recall the main steps of their proof.

We first work away from the hyperplanes $\mathcal D_\lambda$ with $Q(\lambda)>0$ and $c_f^+(-Q(\lambda),\lambda+L)\ne0$.
Let $L_k=-2iv^2\partial/\partial\overline\tau$ be the Maass lowering operator.
Under the dual basis identification $(\overline{\rho_L})^\vee\simeq\rho_L$, the differential identity for the Schwartz forms gives
$$
L_{n+1}\Theta_L^{\mathrm{KM}}(\tau,z)
=i\,dd^c\Theta_L^{\mathrm{U},\vee}(\tau,z).
$$
Using $g=\xi_{1-n}(f)$ and integrating by parts on $\mathcal F_T$, modularity cancels the paired boundary integrals, while integration in $u$ along the upper boundary extracts the constant Fourier coefficient.
Since the contributions from $f^-$ and from nonzero lattice vectors vanish as $T\to\infty$, the zero vector alone contributes to the boundary limit, yielding
$$
\Lambda_{\mathrm U}(g)(z)
=-i\int_{\mathcal F}^{\mathrm{reg}}
\{f(\tau),dd^c\Theta_L^{\mathrm{U},\vee}(\tau,z)\}\,d\mu(\tau)
+c_f^+(0,0)\varphi_{\mathrm{KM}}(0,z).
$$

To pass $dd^c$ through the regularized integral, we subtract the $z$-independent divergent term $c_f^+(0,0)\log T$ and split the constant Fourier coefficient of the integrand into the finitely many potentially singular lattice terms and the remaining sum.
On compact subsets of our open set, the former decay exponentially, while Gaussian estimates give integrable bounds for the remainder and its derivatives in $z$.
Differentiation in $z$ therefore commutes with the regularized limit, and substituting $\varphi_{\mathrm{KM}}(0,z)=i\Omega(z)$ yields $dd^c\widetilde\Phi_L^{\mathrm U}(z,f)=i\Lambda_{\mathrm U}(g)(z)+c_f^+(0,0)\Omega(z)$.

By Proposition \ref{prop:singularThetaRegularity}, both sides are smooth on $\mathcal D_V\setminus|\widetilde H_{\mathrm U}(f)|$, so the identity extends to this complement by continuity and descends to $X_\Gamma\setminus|H_{\mathrm U}(f)|$.
Since the right-hand side is a smooth closed form on all of $X_\Gamma$, it provides the unique smooth extension of $dd^c\Phi_L^{\mathrm U}(\,\cdot\,,f)$.
\end{proof}
\subsection{The weak converse theorem and the converse theorem for Borcherds products}
\label{subsec:borcherdsConverse}
We now consider when a meromorphic automorphic form whose divisor is a linear combination of Heegner divisors arises as a unitary Borcherds product.
Let
$$
P(\tau)=\sum_{\mu\in A_L}\sum_{\substack{m<0\\m+\mathbb Z=-q_L(\mu)}}c(m,\mu)e^{2\pi i m\tau}\mathbf e_\mu
$$
be a finite Fourier polynomial with $c(m,\mu)=c(m,-\mu)\in\mathbb Z$.
Because $1-n<0$, there exists a unique $f_P\in H^+_{1-n,\overline{\rho_L}}$ with principal part $P$ by \cite[Proposition 3.11]{BFtwolift}. 
By \cite[Proposition 1.10 and Remark 1.14]{Bruinier2002book}, $f_P$ is a real linear combination of Bruinier's Maass--Poincar\'e series and has real Fourier coefficients.

The following precise form of Theorem \ref{thm:introWeakConverse} will be proved in Section \ref{sec:weakConverseProof}.
It is the unitary analogue of the weak converse theorem of \cite[Theorem 4.23]{Bruinier2002book}.
\begin{theorem}
\label{thm:unitaryWeakConverse}
Let $P$ be an integral symmetric principal part, and let $F$ be a nonzero meromorphic modular form for $\Gamma$ of rational weight $r$ with a finite order multiplier system.
Suppose that $\operatorname{div}_{X_\Gamma}(F)=H_{\mathrm{U}}(f_P)$.
Then there is a real constant $C$ such that, on $X_\Gamma\setminus|H_{\mathrm{U}}(f_P)|$, we have
$$
\log\lVert F\rVert=-\frac14\Phi_L^{\mathrm{U}}(\,\cdot\,,f_P)+C.
$$
Moreover, the constant coefficient in the zero component of $f_P^+$ satisfies $c_{f_P}^+(0,0)=2r$.
\end{theorem}

The converse assertion in Theorem \ref{thm:introBorcherdsConverse} has the following precise formulation. 
We say that the unitary Borcherds converse theorem holds for $(L,\Gamma)$ if every pair $(P,F)$ satisfying the hypotheses of Theorem \ref{thm:unitaryWeakConverse}, with $F$ of weight $r$, admits $f\in M^!_{1-n,\overline{\rho_L}}$ with integral principal part and $C\in\mathbb C^\times$ such that
$H_{\mathrm U}(f)=H_{\mathrm U}(f_P)$, $c_f(0,0)=2r$, and $F=C\,\widetilde\Psi_{\mathrm U}(\,\cdot\,,f)$.
Assuming Theorem \ref{thm:unitaryWeakConverse}, we obtain the following criterion.
\begin{proposition}
\label{prop:kmConverseCriterion}
If $\Lambda_{\mathrm{U}}$ is injective, then the unitary Borcherds converse theorem holds for $(L,\Gamma)$.
\end{proposition}

\begin{proof}
Fix a pair $(P,F)$ as above, where $F$ has weight $r$.
Theorem \ref{thm:unitaryWeakConverse} yields $c_{f_P}^+(0,0)=2r$ and shows that $\Phi_L^{\mathrm{U}}(\,\cdot\,,f_P)=-4\log\lVert F\rVert+C$ on $X_{\Gamma}\setminus|H_{\mathrm{U}}(f_P)|$ for some real constant $C$.
Since $F$ has weight $r$, we have $dd^c\log\lVert F\rVert=-r\Omega/2$ on $X_{\Gamma}\setminus|H_{\mathrm{U}}(f_P)|$, and Proposition \ref{prop:thetaLiftCurvature} then implies that $\Lambda_{\mathrm{U}}(\xi_{1-n}(f_P))$ vanishes on this dense open subset.
By continuity, the smooth form $\Lambda_{\mathrm{U}}(\xi_{1-n}(f_P))$ vanishes on all of $X_{\Gamma}$.
The injectivity of $\Lambda_{\mathrm{U}}$ and Proposition \ref{prop:xiKernel} therefore imply that $f_P\in M^!_{1-n,\overline{\rho_L}}$, with $c_{f_P}(0,0)=c_{f_P}^+(0,0)=2r$.

Next, Theorem \ref{thm:unitaryBorcherdsProduct} produces a Borcherds product $\widetilde\Psi_{\mathrm{U}}(\,\cdot\,,f_P)$ with the same weight and divisor as $F$, so the quotient of their lifts to $\mathcal D_{V}$ is a nowhere-vanishing holomorphic function.
Comparing the norm identity above with that for $\lVert\widetilde\Psi_{\mathrm{U}}(\,\cdot\,,f_P)\rVert$ at the end of Subsection \ref{subsec:singularThetaLift} shows that this quotient has constant modulus.
Since $\mathcal D_{V}$ is connected, the quotient is constant, which yields $F=C'\widetilde\Psi_{\mathrm{U}}(\,\cdot\,,f_P)$ for some $C'\in\mathbb C^\times$.
\end{proof}
Let $*$ denote the Hodge star operator of the invariant K\"ahler metric, and let $d^*=-*d*$ be the formal adjoint of $d$.
We denote by $\operatorname{tr}_\Omega$ the adjoint of exterior multiplication by $\Omega$ on differential forms.
On $(1,1)$-forms, this is equal to the usual trace with respect to $\Omega$, and $\operatorname{tr}_\Omega\Omega=n$.
A $(1,1)$-form is primitive precisely when its trace vanishes.
Let $\mathcal H_{(2)}^{1,1}(X_{\Gamma})$ denote the space of square integrable harmonic $(1,1)$-forms for the invariant K\"ahler metric.
\begin{proposition}
\label{prop:kmCohomologyKernel}
Let $\Lambda_{\mathrm{U}}^{\mathrm{dR}}$ be the conjugate linear map
$$
\Lambda_{\mathrm{U}}^{\mathrm{dR}}:S_{n+1,\rho_L}\longrightarrow H_{\mathrm{dR}}^2(X_{\Gamma},\mathbb C),
\qquad g\longmapsto[\Lambda_{\mathrm{U}}(g)]_{\mathrm{dR}}.
$$
If $n\geq3$, then $\Lambda_{\mathrm{U}}(g)\in\mathcal H_{(2)}^{1,1}(X_{\Gamma})$ for every $g\in S_{n+1,\rho_L}$, and $\ker\Lambda_{\mathrm{U}}^{\mathrm{dR}}=\ker\Lambda_{\mathrm{U}}$.
\end{proposition}

\begin{proof}
The trace calculation in the proof of Lemma \ref{lem:singularThetaLaplacian} shows that $\operatorname{tr}_{\Omega}\Lambda_{\mathrm{U}}(g)=0$ for every $g\in S_{n+1,\rho_L}$, so $\Lambda_{\mathrm{U}}(g)$ is primitive.
Since it is a closed $(1,1)$-form, we also have $d^c\Lambda_{\mathrm U}(g)=0$.
The K\"ahler identity $d^*=4\pi[\operatorname{tr}_\Omega,d^c]$ therefore shows that $\Lambda_{\mathrm U}(g)$ is coclosed and hence harmonic.
Moreover, the square integrability follows by applying the doubling argument of \cite[Proposition 4.1]{BruinierFunkeKMInjectivity} to the unitary pair $\mathrm U(2,2)\times\mathrm U(V)$.

Let $X_{\Gamma}^{\mathrm{BB}}$ be the Baily--Borel compactification of $X_{\Gamma}$, whose boundary is a finite set of cusps of complex codimension $n$.
For $n\geq3$, restriction induces an isomorphism $IH^2(X_{\Gamma}^{\mathrm{BB}},\mathbb C)\simeq H_{\mathrm{dR}}^2(X_{\Gamma},\mathbb C)$ by \cite[Section 3]{GoreskyMacPhersonIH2}, where $IH^2$ denotes middle perversity intersection cohomology.
The $L^2$-Hodge theorem and the proof of Zucker's conjecture \cite{LooijengaL2,SaperSternL2} identify the reduced $L^2$-cohomology $\overline H_{(2)}^2(X_\Gamma)$ with both the space of square integrable harmonic $2$-forms and $IH^2(X_{\Gamma}^{\mathrm{BB}},\mathbb C)$.
These identifications send a harmonic form to its ordinary de Rham class, so $[\Lambda_{\mathrm{U}}(g)]_{\mathrm{dR}}=0$ implies $\Lambda_{\mathrm{U}}(g)=0$.
For an orbifold quotient, the same conclusion follows on a finite neat cover by taking invariants.
\end{proof}

For $n\geq3$, injectivity of $\Lambda_{\mathrm{U}}^{\mathrm{dR}}$ therefore implies the unitary Borcherds converse theorem by Proposition \ref{prop:kmConverseCriterion}.
\section{Proof of the weak converse theorem}
\label{sec:weakConverseProof}

Let $F$, $P$, and $r$ be as in Theorem \ref{thm:unitaryWeakConverse}.
Choose a neat finite index subgroup $\Gamma'\subseteq\Gamma$, and let $\pi:X_{\Gamma'}\to X_\Gamma$ be the natural quotient map.
Viewed as a modular form of weight $r$ for $\Gamma'$, $F$ has divisor $\pi^*H_{\mathrm{U}}(f_P)$ and Petersson norm given by pulling back its norm from $X_\Gamma$.
Writing $\Phi'=\pi^*\Phi_L^{\mathrm{U}}(\,\cdot\,,f_P)$ for the corresponding singular theta lift on $X_{\Gamma'}$, any relation $\log\lVert F\rVert+\frac14\Phi'=C$ on $X_{\Gamma'}\setminus\pi^{-1}\vert H_{\mathrm{U}}(f_P)\vert$ with $C\in\mathbb R$ descends to $\log\lVert F\rVert+\frac14\Phi_L^{\mathrm{U}}(\,\cdot\,,f_P)=C$ on $X_\Gamma\setminus\vert H_{\mathrm{U}}(f_P)\vert$ by the surjectivity of $\pi$.

Next, for a sufficiently divisible positive integer $N$, the modular form $F^N$ has integral weight $Nr$ and trivial multiplier system.
Since $Nf_P$ has principal part $NP$, we have $f_{NP}=Nf_P$ by the uniqueness of the harmonic weak Maass form.
By the linearity of the theta lift and the identity $\log\lVert F^N\rVert=N\log\lVert F\rVert$, establishing both conclusions of Theorem \ref{thm:unitaryWeakConverse} for $(F^N,NP,Nr)$ immediately yields those for $(F,P,r)$ upon dividing by $N$.
These reductions allow us to assume throughout the proof that $\Gamma$ is neat, $r\in\mathbb Z$, and the multiplier system of $F$ is trivial.

Let $d\mu_{\mathrm{U}}=\Omega^n/n!$ be the volume measure of the invariant K\"ahler metric associated with the form $\Omega$ from Section \ref{sec:borcherdsConverse}.
We normalize the Laplacian on smooth functions by $\Delta_{\mathrm{U}} f=\operatorname{tr}_\Omega(dd^c f)$, so that $d^*df=-4\pi\Delta_{\mathrm U}f$.

On $X_\Gamma\setminus\vert H_{\mathrm{U}}(f_P)\vert$, define $u_F=\log\lVert F\rVert+\frac14\Phi_L^{\mathrm{U}}(\,\cdot\,,f_P)$.
We will show that $u_F$ extends smoothly to $X_\Gamma$ with constant Laplacian and belongs to $L^{p}(X_\Gamma,d\mu_{\mathrm{U}})$ for some $p>1$, which will force $u_F$ to be constant.
\begin{lemma}
\label{lem:singularThetaLaplacian}
For $n\geq2$ and $f\in H^+_{1-n,\overline{\rho_L}}$, one has
$$
\Delta_{\mathrm U}\Phi_L^{\mathrm U}(z,f)=n\,c_f^+(0,0)
$$
on $X_\Gamma\setminus|H_{\mathrm U}(f)|$.
\end{lemma}

\begin{proof}
We first prove that $\operatorname{tr}_\Omega\Lambda_{\mathrm{U}}(g)=0$ for every $g\in S_{n+1,\rho_L}$.

For a fixed $z\in\mathcal D_V$, define
$$
G_z(\tau)=\sum_{\mu\in A_L}
\overline{\theta_{L,\mu}^{\mathrm U}(\tau,z)}\,\mathbf e_\mu,
$$
where $\tau=u+iv\in\mathbb H$.
Then $G_z$ has weight $n-1$ and type $\rho_L$.
Let $e_1,\ldots,e_{n+1}$ be an $h$-orthogonal basis of $V_{\mathbb C}$ with Gram matrix $\operatorname{diag}(1,\ldots,1,-1)$ and $z=\mathbb C e_{n+1}$.
Choose local holomorphic coordinates $w=(w_1,\ldots,w_n)$ centered at $z$ in which the corresponding negative line is generated by $s(w)=e_{n+1}+\sum_{j=1}^n w_je_j$.
Since $-h(s(w),s(w))=1-\sum_{j=1}^n|w_j|^2$, we have
$$
\Omega(z)=-\frac{i}{2\pi}\left.\partial\overline\partial \log\left(1-\sum_{j=1}^n|w_j|^2\right)\right|_{w=0}=\frac{i}{2\pi}\sum_{j=1}^n dw_j\wedge d\overline w_j.
$$
For $x=\sum_{j=1}^na_je_j+be_{n+1} \in V_{\mathbb C}$, set $E_x(\tau,z)=\exp\!\left(2\pi i\tau Q(x)+4\pi v h(x_z,x_z)\right)$.
Then $\varphi_0(x,\tau,z)=vE_x(\tau,z)$, and direct differentiation shows that
$$
2i\partial_\tau E_x=-4\pi\sum_{j=1}^n|a_j|^2E_x,
\qquad
R_{n-1}(vE_x)=\left(n-4\pi v\sum_{j=1}^n|a_j|^2\right)E_x,
$$
where $R_k=2i\partial_\tau+k/v$ is the Maass raising operator of weight $k$.

With our normalization, the explicit formula in \cite[Section 3.2]{FunkeHofmannGreen} identifies the coefficient of $dw_j\wedge d\overline w_j$ in $\varphi_{\mathrm{KM}}(x,\tau,z)$ as $(2v|a_j|^2-1/(2\pi))E_x$.
Since the off-diagonal terms have zero trace, we deduce that
$$
\operatorname{tr}_\Omega\varphi_{\mathrm{KM}}(x,\tau,z)=-2\pi i\sum_{j=1}^n\left(2v|a_j|^2-\frac1{2\pi}\right)E_x =iR_{n-1}(vE_x).
$$
Summing over $L^\vee$, we obtain
$
\operatorname{tr}_\Omega\Theta_L^{\mathrm{KM}}(\tau,z)=iR_{n-1}G_z(\tau).
$

Since $g$ is holomorphic, $\partial_\tau\langle G_z,g\rangle_{A_L}=\langle\partial_\tau G_z,g\rangle_{A_L}$.
We have
$$
2i\partial_\tau\!\left(\langle G_z,g\rangle_{A_L}v^{n-1}\right)=\left\langle 2i\partial_\tau G_z+\frac{n-1}{v}G_z,g \right\rangle_{A_L}v^{n-1}=\langle R_{n-1}G_z,g\rangle_{A_L}v^{n-1}.
$$
Using $d\tau\wedge d\overline\tau=-2i\,v^2d\mu$, we obtain $\langle R_{n-1}G_z,g\rangle_{A_L}v^{n+1}\,d\mu=-d\!\left(\langle G_z,g\rangle_{A_L}v^{n-1}\,d\overline\tau \right)$.
Since $G_z$ and $g$ have weights $n-1$ and $n+1$, respectively, and both transform with the unitary representation $\rho_L$, the one-form $\langle G_z,g\rangle_{A_L}v^{n-1}d\overline\tau$ is $\operatorname{SL}_2(\mathbb Z)$-invariant.
The identified boundary pieces have opposite boundary orientations, so the invariance of this one-form implies that their integrals cancel.
Stokes' theorem leads to
$$
\int_{\mathcal F_T} \langle R_{n-1}G_z,g\rangle_{A_L}v^{n+1}\,d\mu =T^{n-1}\int_{-1/2}^{1/2} \langle G_z(u+iT),g(u+iT)\rangle_{A_L}\,du.
$$

For fixed $z$, Gaussian decay ensures $G_z=O(v)$ and $R_{n-1}G_z=O(1)$ as $v\to\infty$, uniformly in $u$.
Combined with the exponential decay of $g$, these bounds show that the integral over $\mathcal F$ converges absolutely and that the boundary term tends to zero as $T\to\infty$.
Passing to the limit $T\to\infty$, we find that
$$
\int_{\mathcal F} \langle R_{n-1}G_z,g\rangle_{A_L}v^{n+1}\,d\mu=0.
$$
Hence
$$
\operatorname{tr}_\Omega\Lambda_{\mathrm{U}}(g)(z)
=\int_{\mathcal F}\left\langle \operatorname{tr}_\Omega\Theta_L^{\mathrm{KM}}(\tau,z),g(\tau) \right\rangle_{A_L}v^{n+1}\,d\mu
=i\int_{\mathcal F} \langle R_{n-1}G_z(\tau),g(\tau)\rangle_{A_L}v^{n+1}\,d\mu
=0.
$$

As $z$ is arbitrary, $\operatorname{tr}_\Omega\Lambda_{\mathrm{U}}(g)$ vanishes identically on $\mathcal D_{V}$ and hence on $X_\Gamma$.
Taking the $\Omega$-trace of the curvature identity in Proposition \ref{prop:thetaLiftCurvature} for $g=\xi_{1-n}(f)$ and using $\operatorname{tr}_\Omega\Omega=n$ then completes the proof.
\end{proof}
Following \cite{BruinierHowardYang2015}, we recall the Siegel domain coordinates near the cusps.
Fix an isotropic $K$-line $I=K\ell$ with $0\ne\ell\in I\cap L$ and choose an isotropic vector $\widetilde\ell\in V$ satisfying $h(\ell,\widetilde\ell)=1$.
Put $V_0=(K\ell\oplus K\widetilde\ell)^\perp$ and write $h_0=h|_{V_0}$.
Each negative line in $\mathcal D_V$ has a unique generator $s_\ell$ satisfying $h(s_\ell,\ell)=1$.
This generator has the form $s_\ell(\zeta,\tau)=\zeta+\delta_K\tau\ell+\widetilde\ell$, with $\zeta\in V_{0,\mathbb C}$ and $\tau\in\mathbb C$.
Put $N(\zeta,\tau):=-h(s_\ell,s_\ell)=2\sqrt{|D_K|}\operatorname{Im}\tau-h_0(\zeta,\zeta)$.
The condition that the line is negative is equivalent to $N(\zeta,\tau)>0$.
Consequently, the map $(\zeta,\tau)\mapsto[s_\ell(\zeta,\tau)]$ identifies $\mathcal D_V$ with the Siegel domain $\mathcal H_{\ell,\widetilde\ell}=\left\{(\zeta,\tau)\in V_{0,\mathbb C}\times\mathbb H:N(\zeta,\tau)>0\right\}$.

The stabilizer $\mathrm U(V)_\ell$ of $\ell$ acts on this Siegel domain.
Let $\Gamma_\ell=\Gamma\cap\mathrm U(V)_\ell$.
Put $U_\varepsilon(\ell)=\{(\zeta,\tau)\in\mathcal H_{\ell,\widetilde\ell}:N(\zeta,\tau)>1/\varepsilon\}$.
For sufficiently small $\varepsilon>0$, the natural map $\Gamma_\ell\backslash U_\varepsilon(\ell)\longrightarrow X_\Gamma$ is an open embedding.
The normalized generator $s_\ell$ is $\Gamma_\ell$-equivariant and therefore descends to a nowhere vanishing holomorphic section of $\mathcal L$ on this cusp neighborhood, with $\lVert s_\ell\rVert^2=N(\zeta,\tau)$.

The center of $\mathrm U(V)_\ell$ consists of the translations $T_a(\lambda)=\lambda+a\delta_K h(\lambda,\ell)\ell$ for $a\in\mathbb Q$ and $\lambda\in V$.
In Siegel domain coordinates, these act by
$T_a(\zeta,\tau)=(\zeta,\tau+a)$.
The center of $\Gamma_\ell$ is of the form $\Gamma_{\ell,T}=\{T_a:a\in w_I\mathbb Z\}$ for a unique $w_I\in\mathbb Q_{>0}$, called the width of the cusp $I$.

Put $q=\exp(2\pi i\tau/w_I)$.
The map $(\zeta,\tau)\mapsto(q,\zeta)$ induces an isomorphism from $\Gamma_{\ell,T}\backslash U_\varepsilon(\ell)$ to a punctured disc bundle
$$
V_\varepsilon(\ell)=\left\{(q,\zeta)\in\mathbb C\times V_{0,\mathbb C}:0<|q|<\exp\left(-\frac{\pi}{w_I\sqrt{|D_K|}}\left(h_0(\zeta,\zeta)+\frac1\varepsilon\right)\right)\right\}
$$
over $V_{0,\mathbb C}$.
Adding the origin in each fibre gives the disc bundle $\widetilde V_\varepsilon(\ell)$.
The action of $\Gamma_\ell/\Gamma_{\ell,T}$ on
$V_\varepsilon(\ell)$ extends to
$\widetilde V_\varepsilon(\ell)$ and preserves
the boundary divisor $q=0$.
This action is free since $\Gamma$ is neat.
The inclusion $V_\varepsilon(\ell)\subset
\widetilde V_\varepsilon(\ell)$ induces an open embedding $\Gamma_\ell\backslash U_\varepsilon(\ell)\longrightarrow(\Gamma_\ell/\Gamma_{\ell,T})\backslash\widetilde V_\varepsilon(\ell)$.
Gluing the quotient on the right to $X_\Gamma$ along
$\Gamma_\ell\backslash U_\varepsilon(\ell)$ gives a partial
compactification.
Performing this construction at all cusps gives the smooth
toroidal compactification $X_\Gamma^{\mathrm{tor}}$.
Its natural proper morphism to $X_\Gamma^{\mathrm{BB}}$
contracts each boundary component to the corresponding cusp.

Near a boundary point over $I$, choose a polydisc neighborhood $U=\{|q|<\varepsilon\}\times B$ in the coordinates $(q,\zeta)$, where $B\subset V_{0,\mathbb C}$ is a bounded polydisc and $0<\varepsilon<e^{-1}$.
Then $U\cap X_\Gamma=U\setminus\{q=0\}$, and on this set
$$
N(\zeta,\tau)=\frac{w_I\sqrt{|D_K|}}{\pi}(-\log|q|)-h_0(\zeta,\zeta).
$$
Since $h_0(\zeta,\zeta)$ is nonnegative and bounded on $B$, shrinking $\varepsilon$ if necessary yields
$$
\frac{w_I\sqrt{|D_K|}}{2\pi}(-\log|q|)\leq N(\zeta,\tau)\leq\frac{w_I\sqrt{|D_K|}}{\pi}(-\log|q|).
$$

A direct calculation gives $d\mu_{\mathrm U}=C\,N(\zeta,\tau)^{-(n+1)}\,dA(\tau)\,dV(\zeta)$ for a constant $C>0$, where $dA$ is Euclidean area measure on $\mathbb C$, and $dV$ is Euclidean volume measure in fixed complex coordinates on $V_{0,\mathbb C}$.
Combining $dA(\tau)=w_I^2(4\pi^2|q|^2)^{-1}dA(q)$ with the preceding bounds on $N$, we obtain
\begin{equation}
\label{eq:toroidalVolumeBound}
d\mu_{\mathrm U}\leq C_2\frac{dA(q)\,dV(\zeta)}{|q|^2(-\log|q|)^{n+1}}
\end{equation}
for some constant $C_2>0$.
\begin{proposition}
\label{prop:weakConverseRegularity}
For $n\geq2$, the function $u_F$ extends smoothly to a real-valued function on $X_\Gamma$. 
This extension belongs to $L^p(X_\Gamma,d\mu_{\mathrm U})$ for some $p>1$.
\end{proposition}

\begin{proof}
Since $f_P$ has real Fourier coefficients, the substitution $\tau\mapsto-\overline\tau$ replaces the integrand in the regularized theta integral by its complex conjugate.
The symmetry of the truncated fundamental domains under this reflection implies that the regularized lift $\Phi_L^{\mathrm U}(\,\cdot\,,f_P)$ and hence $u_F$ is real-valued on $X_\Gamma\setminus|H_{\mathrm U}(f_P)|$.
Given $x\in\vert H_{\mathrm{U}}(f_P)\vert$, choose a nowhere-vanishing holomorphic section $s$ of $\mathcal L$ on a small open neighborhood $W\subset X_\Gamma$ of $x$ and write $F=F_s s^{\otimes r}$, so that $\log\lVert F\rVert=\log\vert F_s\vert+\frac r2\log\lVert s\rVert^2$.
Because $\operatorname{div}(F)=H_{\mathrm{U}}(f_P)$, Proposition \ref{prop:singularThetaRegularity} shows that the combination $\Phi_L^{\mathrm{U}}(\,\cdot\,,f_P)+2\log|F_s|^2$ extends smoothly from $W\setminus|H_{\mathrm{U}}(f_P)|$ to all of $W$.
Writing
$$
u_F =\frac14\bigl(\Phi_L^{\mathrm{U}}(\,\cdot\,,f_P)+2\log|F_s|^2\bigr) +\frac r2\log\lVert s\rVert^2
$$
on $W\setminus|H_{\mathrm{U}}(f_P)|$, we see that $u_F$ extends smoothly to $W$.
As $x$ varies, these local extensions agree with $u_F$ on the dense open set $W\setminus|H_{\mathrm{U}}(f_P)|$, so by continuity they agree on overlaps and glue to a smooth real-valued function on $X_\Gamma$.

The local decomposition of $\log\lVert F\rVert$ shows that it is locally $L^p$-integrable on $X_\Gamma$ for every finite $p\geq1$.
To establish the integrability of $\log\lVert F\rVert$ near the cusps, fix a boundary point of $X_\Gamma^{\mathrm{tor}}$ over a cusp $I=K\ell$ and let $U=\{|q|<\varepsilon\}\times B$ be a polydisc neighborhood as in the preceding construction, with $0<\varepsilon<e^{-1}$.
By \cite[Lemma 10.6]{BailyBorel}, there exists a holomorphic modular form $G$ of positive integral weight $k$ and trivial multiplier that does not vanish at $I$.
Writing $F=F_{s_\ell}s_\ell^{\otimes r}$ and $G=G_{s_\ell}s_\ell^{\otimes k}$ on $U\cap X_\Gamma$, the Fourier--Jacobi expansion of $G_{s_\ell}$ at $I$ has no terms of negative index, and its coefficient of index zero is a nonzero constant.
Thus $G_{s_\ell}$ extends holomorphically across $q=0$ and is nowhere vanishing after shrinking $U$.

The meromorphic function $F^k/G^r$ on $X_\Gamma$ extends meromorphically to $X_\Gamma^{\mathrm{BB}}$ by \cite[Corollary 10.5]{BailyBorel}.
Pulling back this extension to $U$ along $X_\Gamma^{\mathrm{tor}}\to X_\Gamma^{\mathrm{BB}}$ and multiplying by the nowhere-vanishing holomorphic function $G_{s_\ell}^{\,r}$, we obtain a meromorphic extension of $F_{s_\ell}^{\,k}=(F^k/G^r)G_{s_\ell}^{\,r}$ to $U$.
After shrinking $U$, write $F_{s_\ell}^{\,k}=q^ag/h$, where $a\in\mathbb Z$ and $g,h$ are holomorphic on an open neighborhood of $\overline U$ and not divisible by $q$.
We may also assume that neither $g(q,\cdot)$ nor $h(q,\cdot)$ is identically zero for $|q|\leq\varepsilon$.

Using $\lVert s_\ell\rVert^2=N(\zeta,\tau)$, we obtain
$$
\log\lVert F\rVert=\frac ak\log|q|+\frac1k\bigl(\log|g|-\log|h|\bigr)+\frac r2\log N(\zeta,\tau)
$$
on $(U\cap X_\Gamma)\setminus\{gh=0\}$.
The excluded set has measure zero with respect to $d\mu_{\mathrm U}$.
For every finite $p\geq1$, the maps $q\mapsto\log|g(q,\cdot)|$ and $q\mapsto\log|h(q,\cdot)|$ are continuous into $L^p(B)$.
Indeed, applying the Weierstrass preparation theorem in the first $\zeta$-coordinate $w\in\mathbb C$ after a local linear change of coordinates gives polynomials with bounded roots $\alpha$, so the local estimates follow from the uniform $L^{2p}$-bounds for $\log|w-\alpha|$ on bounded discs.
Fubini's theorem and compactness of $\overline B$ give locally uniform $L^{2p}(B)$-bounds for both logarithms, so their almost everywhere convergence implies the asserted continuity by Vitali's convergence theorem.
By compactness of the closed $q$-disc, we obtain
$$
\sup_{|q|\leq\varepsilon}\int_B\bigl(|\log|g(q,\zeta)||^p+|\log|h(q,\zeta)||^p \bigr)\,dV(\zeta)<\infty.
$$

From the preceding bounds on $N(\zeta,\tau)$, we have $|\log N(\zeta,\tau)| \leq C\bigl(1+\log(-\log|q|)\bigr)$ on $U\cap X_\Gamma$ for some constant $C>0$.
For every finite $p\geq1$, we have
$$
\int_B|\log\lVert F\rVert|^p\,dV(\zeta)\leq C_p\left((-\log|q|)^p+1+\bigl(1+\log(-\log|q|)\bigr)^p\right)\leq 3C_p(-\log|q|)^p
$$
for $0<|q|<\varepsilon$.
Here $C_p>0$ is independent of $q$.
For $1\leq p<n$, the upper bound in \eqref{eq:toroidalVolumeBound} shows that
$$
\begin{aligned}
\int_{U\cap X_\Gamma}|\log\lVert F\rVert|^p\,d\mu_{\mathrm{U}}
&\leq C_2\int_{0<|q|<\varepsilon} \frac{\displaystyle\int_B|\log\lVert F\rVert|^p\,dV(\zeta)}{|q|^2(-\log|q|)^{n+1}}\,dA(q)\\
&\leq C_p'\int_{0<|q|<\varepsilon}\frac{dA(q)}{|q|^2(-\log|q|)^{n+1-p}}\\
&=\frac{2\pi C_p'}{n-p}(-\log\varepsilon)^{p-n}<\infty.
\end{aligned}
$$
Here $C_p'>0$ is a constant, and the last equality follows by integration in polar coordinates.
Together with the local $L^p$-integrability on $X_\Gamma$, these boundary estimates imply that $\log\lVert F\rVert\in L^p(X_\Gamma,d\mu_{\mathrm{U}})$ for $1\leq p<n$.

Finally, since $\Phi_L^{\mathrm{U}}(\,\cdot\,,f_P)\in L^p(X_\Gamma,d\mu_{\mathrm{U}})$ for some $p>1$ by \cite[Proposition 4.5(1)]{BruinierHowardYang2015}, we conclude that $u_F=\log\lVert F\rVert+\frac14\Phi_L^{\mathrm{U}}(\,\cdot\,,f_P)$ also belongs to $L^p(X_\Gamma,d\mu_{\mathrm{U}})$ for some $p>1$.
\end{proof}
\begin{proof}[Proof of Theorem \ref{thm:unitaryWeakConverse}]
By Proposition \ref{prop:weakConverseRegularity}, $u_F$ extends to a smooth real-valued function on $X_\Gamma$ and belongs to $L^p(X_\Gamma,d\mu_{\mathrm{U}})$ for some $p>1$.
On $X_\Gamma\setminus|H_{\mathrm{U}}(f_P)|$, the metric identity $dd^c\log\lVert F\rVert=-\frac r2\Omega$ and Lemma \ref{lem:singularThetaLaplacian} imply
\begin{equation}
\Delta_{\mathrm{U}}u_F=\frac n4\bigl(c_{f_P}^+(0,0)-2r\bigr).
\label{eq:uFLaplacian}
\end{equation}
Since $u_F$ is smooth, this identity holds on all of $X_\Gamma$.

Since $X_\Gamma$ is connected and of finite volume with respect to the complete invariant metric, \cite[Corollary 4.22]{Bruinier2002book} implies that $u_F$ is constant.
From \eqref{eq:uFLaplacian}, we conclude that $c_{f_P}^+(0,0)=2r$.
\end{proof}
\section{Cycle integrals and Dirichlet series}
\label{sec:cycleIntegrals	}

Assume throughout this section that $n>2$.
We compute the integrals of the Kudla--Millson lift over special cycles associated with positive definite $K$-lines and express them as special values of the corresponding Dirichlet series.
Our argument is the unitary analogue of Metzler's method in the orthogonal setting \cite{metzler2025newconversetheoremborcherds,metzler2025symmetricsquaretypelseries}.

\subsection{A geometric integral}
\label{subsec:geometricThetaIntegral}

Retaining the notation of Subsection \ref{subsec:unitaryKMLift} for a Hermitian space $(V',h')$ over $K$ of dimension $d\geq3$ and signature $(d-1,1)$, let $\varphi_f\in\mathcal S(V'(\mathbb A_{K,f}))$ and  $K_f\subset H'(\mathbb A_f)$ be a compact open subgroup fixing $\varphi_f$.

Choose a Haar measure $dh_\infty$ on $H'(\mathbb R)$ whose pushforward to $\mathcal D'$ is $\Omega'^{d-1}$ and $dh_f$ on $H'(\mathbb A_f)$ such that $dh_\infty dh_f$ induces the measure $dh$ fixed in \eqref{eq:adelicThetaIntegral}.
Let $h_1,\ldots,h_r$ and $\Gamma_j$ be the representatives and arithmetic groups in the component decomposition \eqref{eq:shimuraComponents}.
We use the orbifold integration convention for the full groups $\Gamma_j$, with the normalization induced by $dh_\infty$.
Set
$$
J_j(g,\varphi_f):=\int_{\Gamma_j\backslash\mathcal D'}\theta(g;\varphi_{\mathrm{KM}}\otimes\varphi_f) \wedge\Omega'^{d-2},
$$
where the theta form is pulled back along
$[z]\mapsto[(z,h_j)]$.

The following proposition compares the adelic theta integral with the integrals over the connected components of $X_{K_f}$.

\begin{proposition}
\label{prop:adelicComponentIntegrals}
For every $g\in G(\mathbb A)$, the theta integral $I_{V'}(g;\widetilde\varphi_{\mathrm{KM}}\otimes\varphi_f)$ and the integrals $J_j(g,\varphi_f)$ converge absolutely.
Moreover,
\begin{equation}
I_{V'}\bigl(g;\widetilde\varphi_{\mathrm{KM}}\otimes\varphi_f\bigr)=\operatorname{vol}(K_f,dh_f)\int_{X_{K_f}}\theta(g;\varphi_{\mathrm{KM}}\otimes\varphi_f)\wedge\Omega'^{d-2}.
\label{eq:adelicGeometricIntegral}
\end{equation}
The integrals $J_j(g,\varphi_f)$ are independent of $j$, so for every $1\leq j\leq r$ one has
\begin{equation}
I_{V'}\bigl(g;\widetilde\varphi_{\mathrm{KM}}\otimes\varphi_f\bigr)=|C_{K_f}|\operatorname{vol}(K_f,dh_f)J_j(g,\varphi_f).
\label{eq:componentThetaIntegral}
\end{equation}
\end{proposition}

\begin{proof}
Put $\varphi_g=\omega_{V'}(g,1)\bigl(\widetilde\varphi_{\mathrm{KM}}\otimes\varphi_f\bigr)$.
A Gaussian bound at infinity and Theorem \ref{thm:convergentUnitarySiegelWeil} ensure
$$
\int_{[H']}\sum_{x\in V'} |\varphi_g(h^{-1}x)|\,dh<\infty.
$$
Decomposing $[H']$ along $C_{K_f}$ as in the proof of \cite[Proposition 4.17]{Kudla03integrals} and applying \eqref{eq:kmScalarization} with the chosen quotient measures, we obtain
$$
I_{V'}\bigl(g;\widetilde\varphi_{\mathrm{KM}}\otimes\varphi_f\bigr)=\operatorname{vol}(K_f,dh_f)\sum_{j=1}^rJ_j(g,\varphi_f).
$$
This proves \eqref{eq:adelicGeometricIntegral} and the absolute convergence of each integral $J_j(g,\varphi_f)$.

It remains to show that the integrals over the connected components are equal. 
Let $\operatorname{U}(1)$ denote the norm-one torus of $K/\mathbb Q$, and choose a nondegenerate $K$-line $\ell\subset V'$.
With respect to the orthogonal decomposition $V'=\ell\oplus\ell^\perp$, the homomorphism $t\mapsto t\,\operatorname{id}_{\ell}\oplus\operatorname{id}_{\ell^\perp}$ defines a section of $\det:H'\to\operatorname{U}(1)$ over $\mathbb Q$.
By strong approximation for $\operatorname{SU}(V')$ away from infinity, the determinant map induces a bijection $C_{K_f}\xrightarrow{\ \sim\ }\operatorname{U}(1)(\mathbb Q)\backslash\operatorname{U}(1)(\mathbb A_f)/\det K_f$.

Let $\eta$ be a nontrivial character of $C_{K_f}$, viewed through the above identification as a character of $\operatorname{U}(1)(\mathbb A)$ trivial at infinity.
On the component represented by $h_j$, choose representatives $h=(h_\infty,h_jk_f)$ with $k_f\in K_f$.
Since $\eta$ is trivial at infinity and on $\det K_f$, we have $\eta(\det h)=\eta(\det h_j)$.
Quotient integration therefore ensures
$$
I_\eta(g,\varphi_f):=\int_{[H']}\theta_{V'}\bigl(g,h;\widetilde\varphi_{\mathrm{KM}}\otimes\varphi_f\bigr)\eta(\det h)\,dh=\operatorname{vol}(K_f,dh_f)\sum_{j=1}^r\eta(\det h_j)J_j(g,\varphi_f).
$$
By finite Fourier inversion, it suffices to prove that $I_\eta(g,\varphi_f)=0$ for every such $\eta$.

For $x\in V'$, let $H_x\subset H'$ be its stabilizer.
Being an extension of a unitary group by a connected unipotent group, $H_x$ is connected and admits no nontrivial $\mathbb Q$-rational characters, so $[H_x]:=H_x(\mathbb Q)\backslash H_x(\mathbb A)$ has finite invariant volume by \cite[Theorem 5.6(i)]{BorelAdeles}.
Moreover, since $d\geq3$, the orthogonal complement $x^\perp$ contains a nondegenerate $K$-line, which provides a $\mathbb Q$-section of $\det:H_x\to\operatorname{U}(1)$.

Let $\mathcal R$ be a set of representatives for the $H'(\mathbb Q)$-orbits in $V'$.
Unfolding the absolutely convergent theta integral, we obtain
$$
I_\eta(g,\varphi_f)=\sum_{x\in\mathcal R}\int_{H_x(\mathbb A)\backslash H'(\mathbb A)} \varphi_g(h^{-1}x) \left(\int_{[H_x]}\eta(\det(ah))\,da\right)\,d\dot h.
$$
The existence of the $\mathbb Q$-section ensures that $\eta\circ\det$ is nontrivial on $H_x(\mathbb A)$, forcing every inner integral to vanish by the right invariance of the measure on $[H_x]$.
Consequently, $I_\eta(g,\varphi_f)=0$ for every nontrivial $\eta$. Thus $J_j(g,\varphi_f)$ is independent of $j$, and \eqref{eq:componentThetaIntegral} follows from \eqref{eq:adelicGeometricIntegral}.
\end{proof}

For the fixed Hermitian lattice $J\subset V'$, define $K_J$ as in Section \ref{sec:unitary shimura varieties}, with $(V,L)$ replaced by $(V',J)$.
This compact open subgroup is the common stabilizer of the functions $\phi_\mu^J$, $\mu\in A_J$.
Since $\lambda(\widetilde\varphi_{\mathrm{KM}}\otimes\phi_\mu^J)=i\,\Phi_\mu$, Theorem \ref{thm:convergentUnitarySiegelWeil} yields $I_{V'}\bigl(g;\widetilde\varphi_{\mathrm{KM}}\otimes\phi_\mu^J\bigr)=\frac{i}{c_{V'}}E(g,s_0,\Phi_\mu).$

\subsection{Cycle integrals of Kudla--Millson lifts}
\label{subsec:kmCycleIntegrals}
Returning to the notation of Section \ref{sec:unitary shimura varieties}, let $V_1\subset V$ be a positive definite $K$-line with respect to $h$.
Put $V_2=V_1^\perp$, which has signature $(n-1,1)$, and set $L_i=L\cap V_i$ for $i=1,2$.

The lattice $M:=L_1\oplus L_2$ has the same rank as $L$, so it has finite index in $L$.
Due to the trace evenness of $L$, $\mathcal I:=L/M\subset A_M$ is isotropic.
Moreover, $M\subset L\subset L^\vee\subset M^\vee$, $\mathcal I^\perp=L^\vee/M$, and $A_L\simeq\mathcal I^\perp/\mathcal I$.
Write $p:\mathcal I^\perp\to A_L$ for the quotient map.
Define the lifting map $\uparrow_L^M:\mathbb C[A_L]\to\mathbb C[A_M]$ by
$$
\uparrow_L^M\mathbf e_\lambda = \sum_{\substack{\mu\in\mathcal I^\perp\\p(\mu)=\lambda}}\mathbf e_\mu,
\qquad \lambda\in A_L.
$$
By \cite[Proposition 2.1]{MullerIsotropicLift}, this map intertwines the Weil representations $\rho_L$ and $\rho_M$. 
For $g\in S_{n+1,\rho_L}$, set $\widetilde g=\uparrow_L^M g\in S_{n+1,\rho_M}$. 
Proposition 2.6(b)--(c) of \cite{metzler2025newconversetheoremborcherds} then gives
$$
\left\langle\Theta_L^{\mathrm{KM}}(\tau),g(\tau) \right\rangle_{A_L}=\left\langle\Theta_M^{\mathrm{KM}}(\tau),\widetilde g(\tau)\right\rangle_{A_M}.
$$

Write $Q_i=Q|_{V_i}$ for $i=1,2$.
For $\mu_1\in A_{L_1}$, define
$$
\theta_{L_1,\mu_1}(\tau)=\sum_{x_1\in L_1+\mu_1}e\bigl(Q_1(x_1)\tau\bigr),
\qquad
\Theta_{L_1}(\tau)=\sum_{\mu_1\in A_{L_1}}\theta_{L_1,\mu_1}(\tau)\mathbf e_{\mu_1}.
$$
The theta series $\Theta_{L_1}$ is a holomorphic modular form of weight $1$ and type $\rho_{L_1}$.

Put $V_{i,\mathbb C}=V_i\otimes_{K,\sigma}\mathbb C$ for $i=1,2$, and let $\mathcal D_2$ be the domain of negative lines in $V_{2,\mathbb C}$.
We identify $H_2=\operatorname{U}(V_2)$ with the subgroup of $H=\operatorname{U}(V)$ fixing $V_1$ pointwise.
The inclusion $V_2\subset V$ induces an $H_2(\mathbb R)$-equivariant holomorphic embedding $i_{V_1}:\mathcal D_2\hookrightarrow\mathcal D_{V}$ of complex codimension one.

Note that $K_{L_2}$ preserves $\widehat L_2$ and acts trivially on $A_{L_2}$.
Its action on $V_1$ is trivial, so it preserves $\widehat M$ and acts trivially on $A_M\simeq A_{L_1}\oplus A_{L_2}$.
Since $\widehat L$ is the inverse image of $\mathcal I$ under the quotient map $\widehat M^\vee\to A_M$, $K_{L_2}$ preserves $\widehat L$.
The induced action on $A_L\simeq\mathcal I^\perp/\mathcal I$ is also trivial, so $K_{L_2}\subset K_L$.
This implies that $\Gamma_{L_2}^{\mathrm U}$ is contained in the pointwise stabilizer $\Delta_{V_1}:=\{\gamma\in\Gamma_L^{\mathrm U}:\gamma|_{V_1}=\operatorname{id}_{V_1}\}$.
Conversely, any $\gamma\in\Delta_{V_1}$ preserves both $V_2$ and $L_2$.
Since $L_2=L\cap V_2$ is primitive in $L$, duality of the trace pairing ensures that the orthogonal projection $\operatorname{pr}_2:V\to V_2$ maps $L^\vee$ onto $L_2^\vee$.
Every $y_2\in L_2^\vee$ therefore lifts to $y=y_1+y_2\in L^\vee$ with $y_1\in V_1$.
Because $\gamma$ fixes $V_1$ pointwise and acts trivially on $A_L$, we obtain $\gamma y_2-y_2=\gamma y-y\in L\cap V_2=L_2$, which proves that $\Delta_{V_1}=\Gamma_{L_2}^{\mathrm U}$.
Setting $\Gamma_{V_1}:=\Gamma\cap\Delta_{V_1}$, the index $d_\Gamma(V_1):=[\Delta_{V_1}:\Gamma_{V_1}]$ is finite since $\Gamma$ has finite index in $\Gamma_L^{\mathrm U}$.

Define the connected cycle associated with $V_1$ by
$
Z_\Gamma(V_1):=\Gamma_{V_1}\backslash\mathcal D_2.
$
Let $\iota_{V_1}:Z_\Gamma(V_1)\to X_\Gamma$ be the map induced by $i_{V_1}$.
For a smooth $(n-1,n-1)$-form $\alpha$ on $X_\Gamma$ with integrable
pullback, we write $\int_{Z_\Gamma(V_1)}\alpha$ for the integral of
$\iota_{V_1}^*\alpha$ over $Z_\Gamma(V_1)$.

We use the normalizations of Subsection \ref{subsec:geometricThetaIntegral} with $V'=V_2$ and $K_f=K_{L_2}$.
Write $\Omega_2$ for the K\"ahler form on $\mathcal D_2$ and $dh_{2,f}$ for the Haar measure on $H_2(\mathbb A_f)$.
Denote the corresponding invariant measure on $[H_2]:=H_2(\mathbb Q)\backslash H_2(\mathbb A)$ by $dh_2$.
For every integrable $\Delta_{V_1}$-invariant top degree form $\eta$ on $\mathcal D_2$, our orbifold integration convention gives
$$
\int_{\Gamma_{V_1}\backslash\mathcal D_2}\eta
=d_\Gamma(V_1)\int_{\Delta_{V_1}\backslash\mathcal D_2}\eta.
$$
Set
$$
\kappa_{L,\Gamma}(V_1)
:=\frac{i\,d_\Gamma(V_1)}
{c_{V_2}\,|C_{K_{L_2}}|\,\operatorname{vol}(K_{L_2},dh_{2,f})},
$$
where $c_{V_2}$ is the Siegel--Weil constant in Theorem \ref{thm:convergentUnitarySiegelWeil}.

The following computation is the unitary analogue of Metzler's calculation in \cite[Section 6.1]{metzler2025newconversetheoremborcherds}.
\begin{proposition}
\label{prop:cycleIntegralFormula}
Let $g\in S_{n+1,\rho_L}$.
Both integrals below converge absolutely, and
$$
\int_{Z_\Gamma(V_1)}\Lambda_{\mathrm{U}}(g)\wedge\Omega^{n-2}= \kappa_{L,\Gamma}(V_1) \int_{\mathcal F}\left\langle \Theta_{L_1}(\tau)\otimes E_{L_2}\left(\tau,\frac{n-1}{2},n\right), \widetilde g(\tau) \right\rangle_{A_M} v^{n+1}\,\frac{du\,dv}{v^2}.
$$
\end{proposition}

For a positive definite Hermitian space $(W,h_W)$ over $K$, put $\varphi_{0,W}(x)=e^{-2\pi h_W(x,x)}$ for $x\in W\otimes_{K,\sigma}\mathbb C$.
For any Kudla--Millson theta series $\Theta_J^{\mathrm{KM}}(\tau)$ and $\mu\in A_J$, write $\theta_{J,\mu}^{\mathrm{KM}}(\tau)$ for its $\mathbf e_\mu$-coefficient.

\begin{lemma}
\label{lem:kmRestriction}
For $x_i\in V_{i,\mathbb C}$, $i=1,2$, the identity $$i_{V_1}^*\varphi_{\mathrm{KM}, V}(x_1+x_2)=\varphi_{0,V_1}(x_1)\varphi_{\mathrm{KM},V_2}(x_2)$$ holds on $\mathcal D_2$.
In particular, $i_{V_1}^*\Omega=\Omega_2$.
For $\tau\in\mathbb H$ and $\mu_i\in A_{L_i}$, $i=1,2$, one has
\begin{equation}
i_{V_1}^*\theta_{M,\mu_1\oplus\mu_2}^{\mathrm{KM}}(\tau)=\theta_{L_1,\mu_1}(\tau)\,\theta_{L_2,\mu_2}^{\mathrm{KM}}(\tau).
\label{eq:thetaFactorization}
\end{equation}
\end{lemma}

\begin{proof}
With our normalization, \cite[Lemma 2.4]{GarciaSankaran} provides the first identity.
We obtain $i_{V_1}^*\Omega=\Omega_2$ by evaluating it at $x_1=x_2=0$.

For $\tau=u+iv$, the definition of the $\tau$-dependent Schwartz form and the first identity imply
$$
\begin{aligned}
i_{V_1}^*\varphi_{\mathrm{KM},V}(x_1+x_2,\tau)
&=e\bigl(uQ_1(x_1)\bigr)
  \varphi_{0,V_1}(\sqrt v\,x_1)\,
  \varphi_{\mathrm{KM},V_2}(x_2,\tau)\\
&=e\bigl(\tau Q_1(x_1)\bigr)\,
  \varphi_{\mathrm{KM},V_2}(x_2,\tau),
\end{aligned}
$$
where the second equality follows from
$\varphi_{0,V_1}(\sqrt v\,x_1)=\exp(-2\pi vQ_1(x_1))$.
Since $M=L_1\oplus L_2$, summing this identity over $x_i\in L_i+\mu_i$ for $i=1,2$ gives the asserted factorization of the theta series.
\end{proof}
For a complex-valued top degree form $\alpha=a\,\Omega_2^{n-1}$, write $|\alpha|:=|a|\,\Omega_2^{n-1}$.
\begin{lemma}
\label{lem:cycleIntegralConvergence}
For every $\mu_2\in A_{L_2}$, the integral
$$
\int_{Z_\Gamma(V_1)}\theta_{L_2,\mu_2}^{\mathrm{KM}}(\tau)\wedge\Omega_2^{n-2}
$$
converges absolutely and uniformly on $\{\tau\in\mathbb H:\operatorname{Im}\tau\geq v_0\}$ for every \(v_0>0\). 
In particular,
$$
\sup_{\tau\in\mathcal F}\int_{Z_\Gamma(V_1)}\left|\theta_{L_2,\mu_2}^{\mathrm{KM}}(\tau)\wedge\Omega_2^{n-2}\right|<\infty.
$$
Moreover, for every $g\in S_{n+1,\rho_L}$,
$$
\int_{\mathcal F}\int_{Z_\Gamma(V_1)}\left|\left\langle i_{V_1}^*\Theta_M^{\mathrm{KM}}(\tau),\widetilde g(\tau)\right\rangle_{A_M}\wedge\Omega_2^{n-2}\right|v^{n+1}\,\frac{du\,dv}{v^2}<\infty.
$$
\end{lemma}

\begin{proof}
Fix $\mu_2\in A_{L_2}$ and $v_0>0$.
The scalarized Kudla--Millson function is a polynomial times a Gaussian.
Using the explicit action of the Weil representation and absorbing the polynomial factor into a more slowly decaying Gaussian, we obtain
$$
v^{-n/2}\left|\theta\bigl(\jmath(g_\tau),h;\widetilde\varphi_{\mathrm{KM},V_2}\otimes\phi_{\mu_2}^{L_2}\bigr)\right|\leq\theta(1,h;\Psi_{\mu_2})
$$
for all $h\in H_2(\mathbb A)$ and $v\geq v_0$.
Here $\Psi_{\mu_2}$ has finite component $\phi_{\mu_2}^{L_2}$, and its archimedean component is a sufficiently large positive multiple
of the standard Gaussian at the same base point raised to the power
$v_0/2$.

Since $n>2$ and the Witt index of $V_2$ is at most one, the convergence criterion in Theorem \ref{thm:convergentUnitarySiegelWeil} ensures the integrability of $\theta(1,\,\cdot\,;\Psi_{\mu_2})$ over $[H_2]$.
Scalarization, quotient integration, and the finite index identity above therefore yield, uniformly for $v\geq v_0$,
$$
\int_{Z_\Gamma(V_1)}
\left|\theta_{L_2,\mu_2}^{\mathrm{KM}}(\tau)\wedge\Omega_2^{n-2}\right|
\ll d_\Gamma(V_1)\int_{[H_2]}\theta(1,h;\Psi_{\mu_2})\,dh_2<\infty.
$$
Since $v$ is bounded away from zero on $\mathcal F$ and on each compact subset of $\mathbb H$, the preceding pointwise comparison proves the uniform bound and locally uniform convergence.

Due to the boundedness of $\Theta_{L_1}$ on $\mathcal F$, \eqref{eq:thetaFactorization} and the preceding estimate show that the inner integral in the final assertion is $O(\lVert\widetilde g(\tau)\rVert)$, which is integrable over $\mathcal F$ by the exponential decay of $\widetilde g$.
\end{proof}

\begin{proof}[Proof of Proposition \ref{prop:cycleIntegralFormula}]
By Lemma \ref{lem:kmRestriction}, the restricted theta kernel is $\Theta_{L_1}(\tau)\otimes\Theta_{L_2}^{\mathrm{KM}}(\tau)$.
Since $\Theta_{L_1}(\tau)$ is constant along the cycle, it suffices to compute the integral of $\Theta_{L_2}^{\mathrm{KM}}(\tau)$ wedged with $\Omega_2^{n-2}$.

For $\mu_2\in A_{L_2}$, the integral of $\theta_{L_2,\mu_2}^{\mathrm{KM}}(\tau)\wedge\Omega_2^{n-2}$ over $Z_\Gamma(V_1)$ is $d_\Gamma(V_1)$ times its integral over the identity component at level $K_{L_2}$. 
On this component the classical theta form is $v^{-n/2}$ times the adelic theta form. 
Hence \eqref{eq:componentThetaIntegral} gives
$$\int_{Z_\Gamma(V_1)}
  \theta_{L_2,\mu_2}^{\mathrm{KM}}(\tau)\wedge\Omega_2^{n-2}= \frac{d_\Gamma(V_1)v^{-n/2}}{|C_{K_{L_2}}|\operatorname{vol}(K_{L_2},dh_{2,f})}I_{V_2}\bigl(\jmath(g_\tau);\widetilde\varphi_{\mathrm{KM},V_2}\otimes\phi_{\mu_2}^{L_2}\bigr).
$$
The Siegel--Weil section of this Schwartz function is $i\Phi_{\mu_2}$.
Hence Theorem \ref{thm:convergentUnitarySiegelWeil} evaluates the adelic integral as $i\,c_{V_2}^{-1}E(\jmath(g_\tau),s_0,\Phi_{\mu_2})$ at $s_0=(n-1)/2$.
Summing the component identities, applying Proposition \ref{prop:adelicClassicalEisenstein}, and restoring the factor $\Theta_{L_1}(\tau)$ gives
\[
\int_{Z_\Gamma(V_1)}i_{V_1}^*\Theta_M^{\mathrm{KM}}(\tau)\wedge\Omega_2^{n-2}=\kappa_{L,\Gamma}(V_1)\Theta_{L_1}(\tau)\otimes E_{L_2}(\tau,s_0,n).
\]
Using the lifting map identity, $i_{V_1}^*\Omega=\Omega_2$, and Fubini's theorem, we obtain
\[
\begin{aligned}
&\int_{Z_\Gamma(V_1)}\Lambda_{\mathrm{U}}(g)\wedge\Omega^{n-2}\\
&\quad=\int_{Z_\Gamma(V_1)}\int_{\mathcal F} \left\langle i_{V_1}^*\Theta_M^{\mathrm{KM}}(\tau), \widetilde g(\tau)\right\rangle_{A_M} \wedge\Omega_2^{n-2}\,v^{n+1}\frac{du\,dv}{v^2}\\
&\quad=\int_{\mathcal F} \left\langle \int_{Z_\Gamma(V_1)} i_{V_1}^*\Theta_M^{\mathrm{KM}}(\tau)\wedge\Omega_2^{n-2}, \widetilde g(\tau) \right\rangle_{A_M}v^{n+1}\frac{du\,dv}{v^2}\\
&\quad=\kappa_{L,\Gamma}(V_1)\int_{\mathcal F} \left\langle \Theta_{L_1}(\tau)\otimes E_{L_2}(\tau,s_0,n), \widetilde g(\tau)\right\rangle_{A_M} v^{n+1}\frac{du\,dv}{v^2}.
\end{aligned}
\]
The bound in Lemma \ref{lem:cycleIntegralConvergence} justifies Fubini and gives absolute convergence of both integrals by the triangle inequality.
\end{proof}
\subsection{Special values of Dirichlet series}
\label{subsec:rankinSelberg}
Fixing $g\in S_{n+1,\rho_L}$ and a $K$-line $V_1\subset V$ on which $h$ is positive definite, we put $\widetilde g=\uparrow_L^M g\in S_{n+1,\rho_M}$.
For $\operatorname{Re}(t)>1-n/2$, define
$$
R_{L_1}(\widetilde g,t):=\int_{\mathcal F} \left\langle\widetilde g(\tau), \Theta_{L_1}(\tau)\otimes E_{L_2}\left(\tau,\frac{n-1}{2}+\overline t,n\right) \right\rangle_{A_M} v^{n+1}\,\frac{du\,dv}{v^2}.
$$
By \cite[Proposition 5.6]{metzler2025symmetricsquaretypelseries}, the integral converges absolutely and $R_{L_1}(\widetilde g,t)$ is holomorphic for $\operatorname{Re}(t)>1-n/2$.
At $t=0$, Proposition \ref{prop:cycleIntegralFormula} and Hermitian symmetry show that
\begin{equation}
\int_{Z_\Gamma(V_1)}\Lambda_{\mathrm{U}}(g)\wedge\Omega^{n-2}=\kappa_{L,\Gamma}(V_1)\,\overline{R_{L_1}(\widetilde g,0)}.
\label{eq:cycleRankinSelberg}
\end{equation}

To evaluate $R_{L_1}(\widetilde g,t)$, we define, for $\operatorname{Re}(w)>1+n/2$, the Dirichlet series
$$
L_{L_1,0}(\widetilde g,w):=\sum_{\mu_1\in A_{L_1}} \sum_{0\ne x\in L_1+\mu_1}\frac{b_{\widetilde g}\bigl(Q_1(x),\mu_1\oplus0\bigr)}{Q_1(x)^w},
$$
where $b_{\widetilde g}(m,\mu)$ denotes the $(m,\mu)$-th Fourier coefficient of $\widetilde g$.
Here the Fourier indices are admissible because $Q_1(x)+\mathbb Z=q_M((x\bmod L_1)\oplus0)$.

Since the even quadratic lattice $(M_{\mathbb Z},Q)=(L_{1,\mathbb Z},Q_1)\oplus(L_{2,\mathbb Z},Q_2)$ carries the Weil representation $\rho_M$ in Metzler's convention, $L_{L_1,0}(\widetilde g,w)$ coincides with the series $L_{L_1,0,1}(\widetilde g,w)$ in \cite[Definition 5.1]{metzler2025symmetricsquaretypelseries}.
\begin{lemma}
\label{lem:lineDirichletConvergence}
The series $L_{L_1,0}(\widetilde g,w)$ converges absolutely and locally uniformly for $\operatorname{Re}(w)>1+n/2$, and hence is holomorphic in this half-plane.
\end{lemma}

\begin{proof}
We apply \cite[Lemma 5.2]{metzler2025symmetricsquaretypelseries} to the cusp form $\widetilde g$ for the quadratic lattice $(M_{\mathbb Z},Q)$.
Since $L_{1,\mathbb Z}$ has rank $2$ and $M_{\mathbb Z}$ has even rank, the convergence condition for weight $n+1$ is
$$
\operatorname{Re}(w)>\frac{(n+1)+2}{2}-\frac12=1+\frac n2.
$$
\end{proof}

\begin{proposition}
\label{prop:rankinSelbergEvaluation}
For $\operatorname{Re}(t)>1-n/2$, one has
$$
R_{L_1}(\widetilde g,t)=\frac{\Gamma(t+n)}{(4\pi)^{t+n}}L_{L_1,0}(\widetilde g,t+n).
$$
\end{proposition}

\begin{proof}
Up to the normalization of the Eisenstein series, this is \cite[Proposition 5.6]{metzler2025symmetricsquaretypelseries} applied to $\widetilde g$ for the trace lattice $(M_{\mathbb Z},Q)$, with $m_1=2$, $k=n+1$, $k_2=n$, $\eta=0$, and $s=\overline t$ in their notation.
\end{proof}

\begin{corollary}
\label{cor:cycleDirichletValue}
Let $g\in S_{n+1,\rho_L}$ and let $V_1\subset V$ be a $K$-line on which $h$ is positive definite.
With $\widetilde g=\uparrow_L^M g$, one has
\begin{equation}
\int_{Z_\Gamma(V_1)}\Lambda_{\mathrm{U}}(g)\wedge\Omega^{n-2}=\kappa_{L,\Gamma}(V_1)\,\frac{\Gamma(n)}{(4\pi)^n}\overline{L_{L_1,0}(\widetilde g,n)}.
\label{eq:cycleDirichletValue}
\end{equation}
In particular, if $\Lambda_{\mathrm{U}}(g)=0$ as a differential form on $X_{\Gamma}$, then $L_{L_1,0}(\uparrow_L^M g,n)=0$ for every positive definite $K$-line $V_1\subset V$.
\end{corollary}

\begin{proof}
Combining Proposition \ref{prop:rankinSelbergEvaluation} at $t=0$ with \eqref{eq:cycleRankinSelberg} establishes \eqref{eq:cycleDirichletValue}, and the vanishing assertion then follows due to the nonvanishing of $\kappa_{L,\Gamma}(V_1)$.
\end{proof}
\section{Injectivity of the Kudla--Millson lift}
\label{sec:kmInjectivity}

Continue to assume that $n>2$, and retain the fixed finite index subgroup $\Gamma\subseteq\Gamma_L^{\mathrm U}$.
By \eqref{eq:cycleDirichletValue}, vanishing of the Kudla--Millson lift forces the special values of the Dirichlet series attached to all positive-definite $K$-lines to vanish.
Under the $\mathcal{O}_K$ orthogonal dual hyperbolic splitting condition below, these identities imply the vanishing of subseries indexed by arbitrary integral ideals.
Ideal M\"obius inversion gives coprime subseries, and absolute convergence then isolates the norm-one terms, yielding sums of Fourier coefficients over scalar unit orbits.
We also give geometric examples in which the Kudla--Millson lift is injective and a counterexample showing that injectivity can fail without the $\mathcal{O}_K$ hyperbolic splitting condition.
\subsection{Dirichlet series and M\"obius inversion}
\label{subsec:dirichletMobius}

We say that the Hermitian lattice $L$ splits an orthogonal dual hyperbolic plane over $\mathcal{O}_K$ if $L=L_0\perp H_{\mathrm{tr}}$, where $$H_{\mathrm{tr}}=\mathcal O_Ke\oplus\mathcal O_Ke'$$ satisfies $h(e,e)=h(e',e')=0$ and $h(e,e')=\delta_K^{-1}$.
Then the underlying trace lattice of $H_{\mathrm{tr}}$ is self-dual and isometric to $2U$, where $U$ is the even unimodular hyperbolic plane.
Thus $L^\vee=L_0^\vee\perp H_{\mathrm{tr}}$ and $A_L\simeq L_0^\vee/L_0$ as finite quadratic modules.
Since the trace dual lattices $L^\vee$ and $L_0^\vee$ are also $\mathcal O_K$-stable, the $\mathcal O_K$-action on $L^\vee$ descends to $A_L=L^\vee/L$.

For $u=a+b\zeta_K$ and $v=c+d\zeta_K$ with $a,b,c,d\in\mathbb Q$, define
$$
\Delta(u,v)=(\overline u v-\overline v u)/\delta_K =ad-bc\in\mathbb Q.
$$
On $\mathcal O_K$, this is the integral determinant form relative to the ordered $\mathbb Z$-basis $(1,\zeta_K)$.
For $\alpha\in K^\times$, write $N(\alpha)=\alpha\overline\alpha$, and write $N(I)$ for the norm of a nonzero fractional ideal $I$.
A pair $(\mu,m)\in A_L\times\mathbb Q_{>0}$ is admissible if $m+\mathbb Z=q_L(\mu)$.

\begin{lemma}
\label{lem:idealGeneratorsDeterminant}
Let $I$ be a nonzero fractional ideal of $\mathcal O_K$ and $r\in N(I)\mathbb Z\setminus\{0\}$. 
Then $I$ has generators $u,v$ over $\mathcal O_K$ with $\Delta(u,v)=r$.
\end{lemma}

\begin{proof}
Choosing $d\in\mathbb Z_{>0}$ with $dI\subseteq\mathcal O_K$ and replacing $(I,r)$ by $(dI,d^2r)$, we may assume that $I$ is an integral ideal.

Set $t=r/N(I)\in\mathbb Z\setminus\{0\}$ and choose $0\ne a\in tI$.
By the two-generator property of ideals in a Dedekind domain, there exists $0\ne v_0\in I$ such that $I=\mathcal O_Ka+\mathcal O_Kv_0$, and therefore $I=tI+\mathcal O_Kv_0$ since $a\in tI$.
Write $v_0=cv$ with $c\in\mathbb Z_{>0}$ and $v$ primitive in the free $\mathbb Z$-module $I$, so that $I=tI+\mathcal O_Kv$ as well.
Since $v$ is primitive, it extends to a $\mathbb Z$-basis $(u_0,v)$ of $I$ having the same orientation as $(1,\zeta_K)$, for which $\Delta(u_0,v)=[\mathcal O_K:I]=N(I)$.

Now set $u=tu_0$.
Then $\mathcal O_Ku+\mathcal O_Kv=\mathcal O_K tu_0+\mathcal O_K tv+\mathcal O_Kv=tI+\mathcal O_Kv=I$ and $\Delta(u,v)=tN(I)=r$.
\end{proof}

Lemma \ref{lem:idealGeneratorsDeterminant} gives a characterization of this splitting in terms of the trace lattice.

\begin{proposition}
\label{prop:hyperbolicTraceSplitting}
Let $(L,h)$ be a trace even Hermitian lattice of signature $(n,1)$, with $n\geq2$.
Then $L$ splits an orthogonal dual hyperbolic plane over $\mathcal{O}_K$ if and only if $L_{\mathbb Z}$ admits an $\mathcal O_K$-stable orthogonal direct summand $M$ isometric to $2U$.
\end{proposition}

\begin{proof}
The $\mathbb Z$-basis $(e,-\zeta_K e',\zeta_K e,e')$ identifies $(H_{\mathrm{tr}})_{\mathbb Z}$ with $2U$, proving the ``only if'' direction.
Conversely, let $M\subset L$ be an $\mathcal O_K$-stable orthogonal summand with $M_{\mathbb Z}\simeq2U$.
Since trace orthogonality and Hermitian orthogonality to an $\mathcal O_K$-stable sublattice coincide, we have $L=L_1\perp M$, where $L_1$ is positive definite of rank $n-1$.

Choose a hyperbolic $\mathbb Z$-basis $(a_1,b_1,a_2,b_2)$ of $M_{\mathbb Z}$ satisfying $B(a_i,a_j)=B(b_i,b_j)=0$ and $B(a_i,b_j)=\delta_{ij}$.
The only totally $B$-isotropic rational planes in $M\otimes_{\mathbb Z}\mathbb Q$ containing $a_1$ are $\mathbb Qa_1\oplus\mathbb Qa_2$ and $\mathbb Qa_1\oplus\mathbb Qb_2$.
Since $Ka_1$ is such a plane, we may assume after exchanging $a_2$ and $b_2$ if necessary that $Ka_1=\mathbb Qa_1\oplus\mathbb Qa_2$.
Similarly, $Kb_1$ must be either $\mathbb Qb_1\oplus\mathbb Qa_2$ or $\mathbb Qb_1\oplus\mathbb Qb_2$.
Moreover, the first case would force $Ka_1=Kb_1$, contradicting $B(a_1,b_1)=1$.
We therefore have $Kb_1=\mathbb Qb_1\oplus\mathbb Qb_2$. 
With respect to the chosen $\mathbb Z$-basis, one has $M=(M\cap Ka_1)\oplus(M\cap Kb_1)=Ia_1\oplus Jb_1$, where $I=\{\alpha\in K:\alpha a_1\in M\}$ and $J=\{\alpha\in K:\alpha b_1\in M\}$ are nonzero fractional ideals of $\mathcal O_K$.
Put $d=N(I)$.
After rescaling $b_1$ and adjusting $J$ accordingly, we may assume that
$h(a_1,b_1)=(\delta_Kd)^{-1}$.
Since $a_1$ is isotropic, we have $h(\alpha a_1,M)=(\delta_K d)^{-1}\alpha\overline J$ for every $\alpha\in K$, where $\overline J=\{\overline a:a\in J\}$.
Because $M$ is $\mathcal O_K$-stable, we deduce that
$$
\alpha a_1\in M^\vee \quad\Longleftrightarrow\quad h(\alpha a_1,M)\subseteq\delta_K^{-1}\mathcal O_K \quad\Longleftrightarrow\quad \alpha\overline J\subseteq d\mathcal O_K.
$$
As $M$ is trace self-dual, it follows that $I=d\,\overline J^{-1}$, so that $I\overline J=d\mathcal O_K=I\overline I$.
Cancelling the invertible ideal $I$ and taking conjugates, we find that $J=I$.

Now choose a nonzero vector $x_0\in L_1$.
Since $h(x_0,x_0)\in\mathbb Z_{>0}$, by Lemma \ref{lem:idealGeneratorsDeterminant} there exist $u,v\in I$ generating $I$ over $\mathcal O_K$ such that $\Delta(u,v)=d\,h(x_0,x_0)$.
The vector $e=x_0+ua_1+vb_1\in L$ then satisfies
$h(e,e)=h(x_0,x_0)-d^{-1}\Delta(u,v)=0$ and
$h(e,M)=(\delta_Kd)^{-1}(u\overline I+v\overline I)=\delta_K^{-1}\mathcal O_K$.
Choose $f\in M$ with $h(e,f)=\delta_K^{-1}$, and set $e'=f-h(f,f)\zeta_K e\in L$.
Then $h(e',e')=0$ and $h(e,e')=\delta_K^{-1}$, so $H:=\mathcal O_Ke\oplus\mathcal O_Ke'\subset L$ is isometric to $H_{\mathrm{tr}}$.
Finally, because $H_{\mathbb Z}$ is unimodular, the trace pairing $B$ identifies $H$ with $\operatorname{Hom}_{\mathbb Z}(H,\mathbb Z)$.
For every $x\in L$, there is a unique $y\in H$ such that $B(x,z)=B(y,z)$ for all $z\in H$. 
Thus $x-y\in L\cap H^\perp$, which yields $L=H\perp(L\cap H^\perp)$.
\end{proof}
For the remainder of this subsection, assume that
$L=L_0\perp H_{\mathrm{tr}}$.
We next use Lemma \ref{lem:idealGeneratorsDeterminant} to construct positive $K$-lines with prescribed trace dual intersections.
\begin{lemma}
\label{lem:idealLineRepresentatives}
Assume that $L=L_0\perp H_{\mathrm{tr}}$.
Let $(\mu,m)$ be admissible, and let $I\subset\mathcal O_K$ be a nonzero integral ideal.
Then there exists $x\in V$ spanning a positive $K$-line and satisfying $Q(x)=m$, $L^\vee\cap Kx=Ix$, and $\alpha x+L=\alpha\mu$ for every $\alpha\in I$.
\end{lemma}

\begin{proof}
Choose $x_0\in L_0^\vee$ representing $\mu$ with $Q(x_0)>m$, and put $r=Q(x_0)-m\in\mathbb Z_{>0}$.
Since $N(I^{-1})=N(I)^{-1}$, Lemma \ref{lem:idealGeneratorsDeterminant} gives $u,v\in I^{-1}$ generating $I^{-1}$ over $\mathcal O_K$ with $\Delta(u,v)=r$.

Set $x=x_0+ue+ve'$, then it satisfies $Q(x)=Q(x_0)-\Delta(u,v)=m$.
For every $\alpha\in I$, we have $\alpha x\in L^\vee$ and $\alpha x+L=\alpha\mu$.
Conversely, if $\xi\in K$ satisfies $\xi x\in L^\vee$, projection onto $H_{\mathrm{tr}}$ gives $\xi u,\xi v\in\mathcal O_K$.
Since $u,v$ generate $I^{-1}$, it follows that $\xi I^{-1}\subseteq\mathcal O_K$, hence $\xi\in I$.
Thus $L^\vee\cap Kx=Ix$.
\end{proof}
For the positive $K$-lines constructed above, the Dirichlet series of Subsection \ref{subsec:rankinSelberg} can be reindexed by nonzero elements of integral ideals.
Since $q_L(\alpha\mu)=N(\alpha)q_L(\mu)$, the pair $(\alpha\mu,N(\alpha)m)$ is admissible for every admissible $(\mu,m)$ and every $0\ne\alpha\in\mathcal O_K$.

\begin{definition}
\label{def:idealDirichletSeries}
Fix $g\in S_{n+1,\rho_L}$, an admissible pair $(\mu,m)$, and nonzero integral ideals $I,J\subset\mathcal O_K$.
For $s\in\mathbb C$ such that the following sums converge absolutely, define
$$
\begin{aligned}
L_{\mu,m}(g,s)
&=\sum_{\mathclap{0\ne\alpha\in\mathcal O_K}}
\;\frac{b_g\bigl(N(\alpha)m,\alpha\mu\bigr)}{N(\alpha)^s},
\\
L_{\mu,m}^{I}(g,s)
&=\sum_{\mathclap{0\ne\alpha\in I}}
\;\frac{b_g\bigl(N(\alpha)m,\alpha\mu\bigr)}{N(\alpha)^s},
\\
L_{\mu,m}^{(J)}(g,s)
&=\sum_{\mathclap{\substack{0\ne\alpha\in\mathcal O_K\\
(\alpha)+J=\mathcal O_K}}}
\;\frac{b_g\bigl(N(\alpha)m,\alpha\mu\bigr)}{N(\alpha)^s}.
\end{aligned}
$$
\end{definition}

\begin{proposition}
\label{prop:lineIdealSeries}
Fix $g\in S_{n+1,\rho_L}$, an admissible pair $(\mu,m)$, and a nonzero integral ideal $I\subseteq\mathcal O_K$.
Choose $x\in V$ as in Lemma \ref{lem:idealLineRepresentatives}.
Put $V_1=Kx$, $L_1=L\cap V_1$, $L_2=L\cap V_1^\perp$, $M=L_1\oplus L_2$, and $\widetilde g=\uparrow_L^M g$.
For $\operatorname{Re}(s)>1+n/2$, the series $L_{\mu,m}^{I}(g,s)$ converges absolutely and locally uniformly, and satisfies
$$
L_{L_1,0}(\widetilde g,s)=m^{-s}L_{\mu,m}^{I}(g,s).
$$
\end{proposition}
\begin{proof}
By Lemma \ref{lem:lineDirichletConvergence}, $L_{L_1,0}(\widetilde g,s)$ converges absolutely and locally uniformly in the stated half-plane.
For $y\in L_1^\vee$, the natural identification $A_M\simeq A_{L_1}\oplus A_{L_2}$ sends $y+M$ to $(y+L_1)\oplus0$.
Since $(L/M)^\perp=L^\vee/M$, this class lies in $(L/M)^\perp$ if and only if $y\in L^\vee$.
In that case, its image under the natural quotient map $L^\vee/M\to A_L=L^\vee/L$ is $y+L$.
By the definition of the lifting map, the $(y+L_1)\oplus0$-component of $\widetilde g$ is therefore the $y+L$-component of $g$ when $y\in L^\vee$, and is zero otherwise.

Lemma \ref{lem:idealLineRepresentatives} gives $L^\vee\cap V_1=Ix$, together with $Q_1(\alpha x)=N(\alpha)m$ and $\alpha x+L=\alpha\mu$ for every $\alpha\in I$.
Hence
$$
L_{L_1,0}(\widetilde g,s)=\sum_{0\ne y\in L^\vee\cap V_1}\frac{b_g\bigl(Q_1(y),y+L\bigr)}{Q_1(y)^s}=m^{-s}\sum_{0\ne\alpha\in I}\frac{b_g\bigl(N(\alpha)m,\alpha\mu\bigr)}{N(\alpha)^s}=m^{-s}L_{\mu,m}^{I}(g,s).
$$
The termwise comparison also proves the asserted absolute and locally uniform convergence.
\end{proof}
\begin{corollary}
\label{cor:idealSubseriesVanishing}
Assume that $n>2$ and that $L$ splits off $H_{\mathrm{tr}}$.
Let $g\in S_{n+1,\rho_L}$ satisfy $\Lambda_{\mathrm U}(g)=0$.
For every admissible pair $(\mu,m)$ and every nonzero integral ideal $I\subseteq\mathcal O_K$, the series $L_{\mu,m}^{I}(g,n)$ converges absolutely to zero.
\end{corollary}

\begin{proof}
This follows from Proposition \ref{prop:lineIdealSeries} and Corollary \ref{cor:cycleDirichletValue}.
\end{proof}
We now deduce the vanishing of the coprime subseries by M\"obius inversion.

\begin{lemma}
\label{lem:coprimeSubseriesVanishing}
Assume that $n>2$ and that $L$ splits off $H_{\mathrm{tr}}$.
Let $g\in S_{n+1,\rho_L}$ satisfy $\Lambda_{\mathrm U}(g)=0$.
For every admissible pair $(\mu,m)$ and every nonzero integral ideal $J\subseteq\mathcal O_K$, the series $L_{\mu,m}^{(J)}(g,n)$ converges absolutely to zero.
\end{lemma}

\begin{proof}
Taking $I=\mathcal O_K$ in Corollary \ref{cor:idealSubseriesVanishing} gives absolute convergence of $L_{\mu,m}(g,n)$, and hence of its subseries $L_{\mu,m}^{(J)}(g,n)$.

Let $\mu_K$ denote the M\"obius function on the nonzero integral ideals of $\mathcal O_K$, defined by
$$
\mu_K(I)=
\begin{cases}
1, & I=\mathcal O_K,\\
(-1)^r, & I=\mathfrak p_1\cdots\mathfrak p_r,\quad r\geq1,\\
0, & \mathfrak p^2\mid I\text{ for some nonzero prime ideal }\mathfrak p,
\end{cases}
$$
where $\mathfrak p_1,\ldots,\mathfrak p_r$ in the second case are distinct nonzero prime ideals.

For $0\ne\alpha\in\mathcal O_K$, one has
$$
\mathbf1_{(\alpha)+J=\mathcal O_K}=\sum_{\substack{I\mid J\\\alpha\in I}}\mu_K(I).
$$
Substituting this identity and summing over the finitely many divisors of $J$ gives
$$
L_{\mu,m}^{(J)}(g,n)=\sum_{I\mid J}\mu_K(I)L_{\mu,m}^{I}(g,n)=0,
$$
where the last equality follows from Corollary \ref{cor:idealSubseriesVanishing}.
\end{proof}

\subsection{The kernel of the Kudla--Millson lift}
\label{subsec:kmKernel}

For $u\in\mathcal O_K^\times$, define $R_u$ on $\mathbb C[A_L]$ by $R_u\mathbf e_\mu=\mathbf e_{u\mu}$.
Since $N(u)=1$, multiplication by $u$ induces an automorphism of $(A_L,q_L)$, so $R_u$ commutes with $\rho_L$.
It therefore acts on $S_{n+1,\rho_L}$, with $b_{R_ug}(m,\mu)=b_g(m,u^{-1}\mu)$.

The symmetry $b_g(m,-\mu)=b_g(m,\mu)$ implies that $R_{-1}=\operatorname{id}$ on $S_{n+1,\rho_L}$, so the unit action factors through $\overline{\mathcal U}=\mathcal U/\{\pm1\}$.
For $[u]\in\overline{\mathcal U}$, write $R_{[u]}$ for the induced operator.
The averaging operator
$$
P_{\mathcal U}=\frac{1}{|\mathcal U|}\sum_{u\in\mathcal U}R_u
$$
is the projection onto the $\mathcal U$-invariant subspace.

The Kudla--Millson lift is invariant under this action.
Since $u\,\operatorname{id}_V$ fixes $\mathcal D_V$ pointwise and $N(u)=1$, one has $\varphi_{\mathrm{KM}}(ux,\tau,z)=\varphi_{\mathrm{KM}}(x,\tau,z)$.
Reindexing the theta series by multiplication by $u$ therefore gives $R_u\Theta_L^{\mathrm{KM}}(\tau,z)=\Theta_L^{\mathrm{KM}}(\tau,z)$.
As $R_u$ is unitary for the coefficient pairing, this implies $\Lambda_{\mathrm{U}}(R_ug)=\Lambda_{\mathrm{U}}(g)$.
Averaging over $\mathcal U$ yields
$$
\Lambda_{\mathrm{U}}=\Lambda_{\mathrm{U}}\circ P_{\mathcal U}.
$$
In particular, $\ker(P_{\mathcal U})\subseteq\ker(\Lambda_{\mathrm{U}})$.
To prove the reverse inclusion, we adapt Metzler's coprime subseries argument to ideals of $\mathcal O_K$.
\begin{theorem}
\label{thm:unitaryKMKernel}
Suppose that $n>2$ and that $L$ splits an orthogonal dual hyperbolic plane over $\mathcal{O}_K$.
Then the Kudla--Millson lift $\Lambda_{\mathrm{U}}:S_{n+1,\rho_L}\to\mathcal Z^{1,1}(X_{\Gamma})$ satisfies
$$
\ker(\Lambda_{\mathrm{U}})=\ker(P_{\mathcal U}).
$$
\end{theorem}

\begin{proof}
We already know $\ker(P_{\mathcal U})\subseteq\ker(\Lambda_{\mathrm{U}})$.
To prove the reverse inclusion, let $g\in\ker(\Lambda_{\mathrm{U}})$ and suppose that $P_{\mathcal U}g\ne0$.
Choose an admissible pair $(\mu,m)$ with $b_{P_{\mathcal U}g}(m,\mu)\ne0$, and set
$$
A_{\mu,m}(g):=\sum_{u\in\mathcal U}b_g(m,u\mu)=|\mathcal U|\,b_{P_{\mathcal U}g}(m,\mu)\ne0.
$$
Since $L_{\mu,m}(g,n)$ converges absolutely, we may choose $R>1$ such that
$$
\sum_{\substack{\alpha\in\mathcal O_K\setminus\{0\}\\
N(\alpha)>R}}\frac{|b_g(N(\alpha)m,\alpha\mu)|}{N(\alpha)^n}<|A_{\mu,m}(g)|.
$$

Only finitely many $\alpha\in\mathcal O_K$ satisfy $1<N(\alpha)\le R$, so the product
$$
I_R:=\prod_{\substack{\alpha\in\mathcal O_K\setminus\{0\}\\
1<N(\alpha)\le R}}(\alpha)
$$
is a nonzero integral ideal, with the empty product interpreted as $\mathcal O_K$.
For $1<N(\alpha)\le R$, one has $I_R\subseteq(\alpha)$, so $(\alpha)+I_R=(\alpha)\ne\mathcal O_K$ and no such term appears in $L_{\mu,m}^{(I_R)}(g,n)$.
In $L_{\mu,m}^{(I_R)}(g,n)$, the norm-one terms contribute $A_{\mu,m}(g)$, while all remaining terms have norm strictly greater than $R$.
Consequently,
$$
\left|L_{\mu,m}^{(I_R)}(g,n)-A_{\mu,m}(g)\right|
\le\sum_{\substack{\alpha\in\mathcal O_K\setminus\{0\}\\ N(\alpha)>R}}\frac{|b_g(N(\alpha)m,\alpha\mu)|}{N(\alpha)^n}
<|A_{\mu,m}(g)|.
$$
On the other hand, $L_{\mu,m}^{(I_R)}(g,n)=0$ by Lemma \ref{lem:coprimeSubseriesVanishing}, a contradiction.
We conclude that $P_{\mathcal U}g=0$.
\end{proof}
Under the same hypotheses, Theorem \ref{thm:unitaryKMKernel} and Proposition \ref{prop:kmCohomologyKernel} give
$$
\ker(\Lambda_{\mathrm{U}}^{\mathrm{dR}})=\ker(\Lambda_{\mathrm{U}})=\ker(P_{\mathcal U}).
$$

\begin{corollary}
\label{cor:kmInjectivityCriteria}
Assume that $n>2$ and that $L$ splits an orthogonal dual hyperbolic plane over $\mathcal{O}_K$.
Then the lift $\Lambda_{\mathrm{U}}$ is injective on $\operatorname{im}(P_{\mathcal U})=S_{n+1,\rho_L^{\mathrm U}}$. 
Thus the Kudla--Millson lift from vector-valued cusp forms for $\mathrm U(1,1)(\mathbb Z)$ is injective.
It is injective on $S_{n+1,\rho_L}$ if and only if $P_{\mathcal U}=\operatorname{id}$.
In particular, it is injective when $\mathcal O_K^\times=\{\pm1\}$.
\end{corollary}
A sufficient condition for injectivity can be expressed directly in terms of the action of units on $L^\vee/L$.

\begin{corollary}
\label{cor:discriminantUnitCriterion}
Assume that $n>2$ and that $L$ splits an orthogonal dual hyperbolic plane over $\mathcal{O}_K$.
If $(u-1)L^\vee\subseteq L$ for a unit $u\in\mathcal O_K^\times$ such that $\mathcal U=\langle-1,u\rangle$, then $\Lambda_{\mathrm{U}}$ is injective on $S_{n+1,\rho_L}$.
In particular, it suffices that $(i-1)L^\vee\subseteq L$ when $K=\mathbb Q(i)$ and $(\omega_3-1)L^\vee\subseteq L$ when $K=\mathbb Q(\sqrt{-3})$, where $\omega_3$ is a primitive cube root of unity.
\end{corollary}

\begin{proof}
The inclusion $(u-1)L^\vee\subseteq L$ means that $u$ acts trivially on $A_L=L^\vee/L$.
Hence $R_{[u]}$ is the identity on $S_{n+1,\rho_L}$.
Since $[u]$ generates $\overline{\mathcal U}$, we have $P_{\mathcal U}=\operatorname{id}$.
Corollary \ref{cor:kmInjectivityCriteria} gives the asserted injectivity.
The two special cases follow because $[i]$ and $[\omega_3]$ generate $\overline{\mathcal U}$ for the respective fields.
\end{proof}

For $\chi\in\widehat{\overline{\mathcal U}}$, let $S_{n+1,\rho_L}[\chi]$ be the subspace of forms $g$ satisfying $R_{[u]}g=\chi([u])g$ for every $[u]\in\overline{\mathcal U}$.
If $L$ splits an orthogonal dual hyperbolic plane over $\mathcal{O}_K$, the character decomposition gives
$$
\begin{aligned}
S_{n+1,\rho_L}
&=\bigoplus_{\chi\in\widehat{\overline{\mathcal U}}}S_{n+1,\rho_L}[\chi],\\
\ker(\Lambda_{\mathrm U})=\ker(P_{\mathcal U})
&=\bigoplus_{\substack{\chi\in\widehat{\overline{\mathcal U}}\\\chi\ne1}}S_{n+1,\rho_L}[\chi].
\end{aligned}
$$

\begin{corollary}
\label{cor:unitaryBorcherdsConverse}
Assume that $n>2$ and that $L$ splits an orthogonal dual hyperbolic plane over $\mathcal{O}_K$.
If $S_{n+1,\rho_L}[\chi]=0$ for every nontrivial character $\chi\in\widehat{\overline{\mathcal U}}$, then the unitary Borcherds converse theorem holds for $(L,\Gamma)$. 
In particular, it holds when $\mathcal O_K^\times=\{\pm1\}$.
\end{corollary}

To compute $\dim S_{n+1,\rho_L}[\chi]$, let $W_\chi\subseteq\mathbb C[A_L]$ be the subspace on which $R_u$ acts as $\chi([u])$ for every $u\in\mathcal O_K^\times$.
Since the operators $R_u$ commute with $\rho_L$, the subspace $W_\chi$ is stable under this representation.
Write $\rho_\chi=\rho_L|_{W_\chi}$ and $d_\chi=\dim W_\chi$.
By definition, $S_{n+1,\rho_L}[\chi]$ consists of the cusp forms taking values in $W_\chi$, so $S_{n+1,\rho_L}[\chi]=S_{n+1,\rho_\chi}$.
For a unitary operator $U_0$ with eigenvalues $e^{2\pi i\nu_j}$, where $0\leq\nu_j<1$, put $\alpha(U_0)=\sum_j\nu_j$.
Define
$$
\begin{aligned}
a_{1,\chi}
&=\alpha\!\left(e^{\pi i (n+1)/2}\rho_\chi(S)\right),
&
a_{2,\chi}
&=\alpha\!\left(\left(e^{\pi i (n+1)/3}\rho_\chi(ST)\right)^{-1}\right),\\
a_{3,\chi}
&=\alpha\!\left(\rho_\chi(T)\right),
&
a_{4,\chi}
&=\dim\ker\!\left(\rho_\chi(T)-\operatorname{id}_{W_\chi}\right).
\end{aligned}
$$

For $a\in\mathcal O_K$, put $A_L[a]=\{\mu\in A_L:a\mu=0\}$.

\begin{proposition}
Assume that $n>2$.
For every $\chi\in\widehat{\overline{\mathcal U}}$, one has
\begin{equation}
\dim S_{n+1,\rho_L}[\chi]=\frac{n+13}{12}d_\chi-a_{1,\chi}-a_{2,\chi}-a_{3,\chi}-a_{4,\chi}.
\label{eq:unitCharacterDimension}
\end{equation}
For the nontrivial characters in the two exceptional fields, $d_\chi$ is given as follows:
\begin{enumerate}[(i)]
\item
If $K=\mathbb Q(i)$ and $\varepsilon$ is the nontrivial character of $\overline{\mathcal U}=\{[1],[i]\}$, then
$$
d_\varepsilon=\frac{|A_L|+|A_L[2]|-2|A_L[1-i]|}{4}.
$$

\item
If $K=\mathbb Q(\sqrt{-3})$, $\omega_3$ is a primitive cube root of unity, and $\varepsilon$ is a nontrivial character of $\overline{\mathcal U}=\langle[\omega_3]\rangle$, then
$$
d_\varepsilon=d_{\overline{\varepsilon}}=\frac{|A_L|+|A_L[2]|-|A_L[1-\omega_3]|-1}{6}.
$$
\end{enumerate}
\end{proposition}

\begin{proof}
Since $R_{-1}$ acts trivially on $W_\chi$, one has $\rho_\chi(-I)=(-1)^{n+1}\operatorname{id}_{W_\chi}$, and the representation $\rho_\chi$ is unitary with finite image.
By the Riemann--Roch dimension formula in \cite[p.~228]{Borcherds1999GrossKohnenZagier}, we have
$$
\dim M_{n+1,\rho_\chi}=\frac{n+13}{12}d_\chi-a_{1,\chi}-a_{2,\chi}-a_{3,\chi}.
$$
Because the constant term map $M_{n+1,\rho_\chi}\to\ker(\rho_\chi(T)-\operatorname{id}_{W_\chi})$ is surjective with kernel $S_{n+1,\rho_\chi}$, we have $\dim S_{n+1,\rho_\chi}=\dim M_{n+1,\rho_\chi}-a_{4,\chi}$, which establishes \eqref{eq:unitCharacterDimension}.

Next, since $R_u\mathbf e_\mu=\mathbf e_{u\mu}$, the trace of $R_u$ on $\mathbb C[A_L]$ equals the number of fixed points of $u$, namely $|A_L[1-u]|$.
Viewing $\chi$ as the character $u\mapsto\chi([u])$ of the full unit group $\mathcal O_K^\times$, we obtain from the character projector formula that
$$
d_\chi=\frac{1}{|\mathcal O_K^\times|}\sum_{u\in\mathcal O_K^\times}\overline{\chi([u])}\,|A_L[1-u]|.
$$
Evaluating this sum for $K=\mathbb Q(i)$ and $K=\mathbb Q(\sqrt{-3})$ directly proves parts~(i) and~(ii).
\end{proof}

As an application of our injectivity results, we obtain information about the Heegner part of Picard
groups.
Let $\operatorname{Pic}_{\mathrm{Heeg}}(X_L)_{\mathbb Q}$ denote the subspace of $\operatorname{Pic}(X_L)\otimes_{\mathbb Z}\mathbb Q$ spanned by the classes of unitary Heegner divisors.
Arguing as in \cite[Corollary 7.1]{Bruinier14Converse},
after averaging over $\mathcal U$, we obtain the following corollary.
\begin{corollary}
\label{cor:heegnerPicardDimension}
Suppose that $n>2$ and that $L$ splits an orthogonal dual hyperbolic plane over $\mathcal O_K$.
Then
$$
\dim_{\mathbb Q}\operatorname{Pic}_{\mathrm{Heeg}}(X_L)_{\mathbb Q}
=1+\dim_{\mathbb C}S_{n+1,\rho_L}^{\mathcal U}.
$$
In particular, if $P_{\mathcal U}=\operatorname{id}$ on $S_{n+1,\rho_L}$, this dimension equals $1+\dim_{\mathbb C}S_{n+1,\rho_L}$.
\end{corollary}

\subsection{Examples over the Eisenstein integers}
\label{subsec:eisensteinLatticeExamples}
We now describe some moduli spaces associated with ball quotients whose lattices satisfy the hypotheses of Corollary \ref{cor:discriminantUnitCriterion}.
Let $K=\mathbb Q(\sqrt{-3})$ and write $\mathcal O_K=\mathbb Z[\omega_3]$, where $\omega_3$ is a primitive cube root of unity.
Let $\Lambda_4$ be the positive definite Eisenstein lattice obtained by multiplying the Hermitian form of a rank-four block in \cite[Theorem 1.7]{AllcockCarlsonToledoCubicThreefolds} by $1/3$, so that its trace lattice is $E_8$ (see also \cite[\S5]{Allcock2000Leech}).
Let $\Lambda_1$ be the rank-one Eisenstein lattice whose trace lattice is the positive definite root lattice $A_2$.

\begin{enumerate}[(i)]
\item
By \cite[Theorem 6.1]{AllcockCarlsonToledoCubicThreefolds}, the period map identifies the moduli space of smooth cubic threefolds with the complement of the nodal and chordal Heegner divisors in the ball quotient $\mathrm{U}(L_{\mathrm{cub}})\backslash\mathbb B^{10}$ (see also \cite[Theorem 3.1]{LooijengaSwierstraCubic}).
Multiplying the Hermitian form in \cite[Theorem 1.7]{AllcockCarlsonToledoCubicThreefolds} by $1/3$ to match our normalization of $h$ leaves the ball quotient unchanged.
In this normalization, the Hermitian lattice is
$$
L_{\mathrm{cub}}=\Lambda_1\perp\Lambda_4^{\oplus2}\perp H_{\mathrm{tr}}.
$$
Since $A_{L_{\mathrm{cub}}}\simeq\mathbb Z/3\mathbb Z$ and the scalar $-1$ induces its nontrivial isometry, $\mathrm{U}(L_{\mathrm{cub}})$ and $\Gamma_{L_{\mathrm{cub}}}^{\mathrm U}$ have the same projective image and hence the same coarse ball quotient.

\item
By \cite{Mostow1986GeneralizedPicard,Allcock2000Leech}, the Deligne--Mostow period map for twelve equal weights $1/6$ identifies the moduli space of unordered configurations of twelve distinct points on $\mathbb P^1$ up to projective equivalence with the complement in $\mathrm{U}(L_{12})\backslash\mathbb B^9$ of the Heegner divisor corresponding to configurations with coincident points.
Taking the double cover of $\mathbb P^1$ branched at these twelve points identifies this space with the moduli space of smooth hyperelliptic curves of genus five.
In our normalization, we may take
$$
L_{12}=\Lambda_4^{\oplus2}\perp H_{\mathrm{tr}}.
$$
Since $A_{L_{12}}=0$, we have $\mathrm{U}(L_{12})=\Gamma_{L_{12}}^{\mathrm U}$.
For this lattice, the source of the Kudla--Millson lift is $S_{10,\rho_{L_{12}}}=S_{10}(\operatorname{SL}_2(\mathbb Z))=0$.
\end{enumerate}

In both cases, $\omega_3$ acts trivially on $A_L$ because the summands $\Lambda_4$ and $H_{\mathrm{tr}}$ are trace self-dual, while the discriminant group of $\Lambda_1$ has order $3$.
Corollaries \ref{cor:discriminantUnitCriterion} and \ref{cor:unitaryBorcherdsConverse} therefore imply that the Kudla--Millson lift is injective and that the unitary Borcherds converse theorem holds on each corresponding discriminant kernel quotient.
\subsection{A counterexample without splitting an orthogonal dual hyperbolic plane over $\mathcal{O}_K$}
\label{subsec:splittingCounterexample}

We construct a Hermitian lattice whose Kudla--Millson lift is not injective.
In fact, the lattice splits the standard Hermitian hyperbolic plane with Gram matrix $\left(\begin{smallmatrix}0&1\\1&0\end{smallmatrix}\right)$.
Its trace lattice splits off $2U$ abstractly and also admits an $\mathcal O_K$-stable orthogonal direct summand isometric to $U\perp U(7)$.

Let $K=\mathbb Q(\sqrt{-7})$, put $\delta_K=\sqrt{-7}$ and $\omega=(1+\delta_K)/2$, and write $\mathcal O_K=\mathbb Z[\omega]$.
Let $H_{h}=\mathcal O_Ke\oplus\mathcal O_Ke'$ with $h(e,e)=h(e',e')=0$ and $h(e,e')=h(e',e)=1$.
Let $u_1,\ldots,u_4$ be an orthonormal basis of a positive definite Hermitian space over $K$, and set
$$
L=H_h\perp \mathcal O_Ku_1\perp\cdots\perp\mathcal O_Ku_4.
$$
Then $(L,h)$ has signature $(5,1)$ and is self-dual with respect to $h$.
We have $L^\vee=\delta_K^{-1}L$.

The $\mathbb Z$-basis $(e,\omega e',(\omega-4)e,(2\omega-1)e')$ identifies $(H_h)_{\mathbb Z}$ with $U\perp U(7)$, which forms an $\mathcal O_K$-stable orthogonal summand of $L_{\mathbb Z}$.
In the orthogonal complement of the first $U$, set $x_0=(\omega-4)e-(2\omega-1)e'+u_1+\omega(u_2+u_3+u_4)$ and $y_0=\omega u_1-2x_0$.
Since $Q(x_0)=Q(y_0)=0$ and $B(x_0,y_0)=1$, these vectors span a second copy of $U$.
By unimodularity, the two hyperbolic planes split off orthogonally, so that $L_{\mathbb Z}\simeq 2U\perp N_1$ for an even positive definite lattice $N_1$ of rank $8$.

Nevertheless, $L$ contains no orthogonal dual hyperbolic plane over $\mathcal{O}_K$.
Indeed, if $M\subseteq L$ were such a plane, then $h(M,M)\subseteq\mathcal O_K$ would imply $\delta_K^{-1}M\subseteq M^\vee=M$.
This contradicts $[\delta_K^{-1}M:M]=7^2$.

\begin{proposition}
\label{prop:noninjectiveExample}
For the lattice $L$ defined above, there exists a nonzero $G\in S_{6,\rho_L}$ such that $\Lambda_{\mathrm U}(G)=0$ on $X_\Gamma$ for every finite index subgroup $\Gamma\subseteq\Gamma_L^{\mathrm U}$.
\end{proposition}

\begin{proof}
We first prove that $\mathbf e_0$ is fixed by $\rho_L(\Gamma_0(7))$.
Multiplication by $\delta_K$ identifies $A_L$ with $L/\delta_KL\simeq\mathbb F_7^6$.
The identity
$$
q_L(\delta_K^{-1}x+L)=\frac{h(x,x)}7+\mathbb Z,
\qquad x\in L,
$$
shows that $(A_L,q_L)$ has level $7$.
The Weil representation attached to $(A_L,-q_L)$ in \cite[Section 2]{Scheithauer2015Constructions} is $\rho_L$ in our convention.
For $\gamma=\left(\begin{smallmatrix}a&b\\c&d\end{smallmatrix}\right) \in\Gamma_0(7)$, we have
$$
\rho_L(\gamma)\mathbf e_0=\left(\frac{a}{7^6}\right)\mathbf e_0=\mathbf e_0.
$$

Set $f_7(\tau)=\eta(\tau)^6\eta(7\tau)^6$, where $\eta$ is the Dedekind eta function.
By \cite[Theorems 1.64 and 1.65]{OnoWebModularity}, we have $f_7\in S_6(\Gamma_0(7))$ with trivial character.
For $W_7=\left(\begin{smallmatrix}0&-1\\7&0\end{smallmatrix}\right)$, we use the normalized slash operator
$$
(F|_6W_7)(\tau)=7^{-3}\tau^{-6}F\!\left(-\frac1{7\tau}\right).
$$
By the transformation law for $\eta$, one has $f_7|_6W_7=-f_7$.

By \cite[Theorem 3.1(i)]{Scheithauer2015Constructions}, the function
$$
G(\tau)= \sum_{\gamma\in\Gamma_0(7)\backslash\operatorname{SL}_2(\mathbb Z)}(f_7|_6\gamma)(\tau)\,\rho_L(\gamma^{-1})\mathbf e_0
$$
is a modular form of weight $6$ and type $\rho_L$.
Since $f_7$ vanishes at every cusp of $\Gamma_0(7)$, each summand decays exponentially at $i\infty$, so $G\in S_{6,\rho_L}$.
For $r=0,\ldots,6$, the classes $\mu_r=\delta_K^{-1}(e+r\omega e')+L$ are nonzero and satisfy $q_L(\mu_r)=r/7+\mathbb Z$.
The component formula in \cite[Theorem 5.7]{Scheithauer2009Weil} therefore shows that the map from $S_6(\Gamma_0(7))$ to $S_{6,\rho_L}$ defined by the preceding coset sum is injective.
Hence $G\ne0$, since $f_7\ne0$.

It remains to prove that $\Lambda_{\mathrm{U}}(G)=0$.
Fix $z\in\mathcal D_{V}$, and write $\tau=u+iv$.
Let $\Theta_{L,0}^{\mathrm{KM}}$ denote the $\mathbf e_0$-component of $\Theta_L^{\mathrm{KM}}$.
Substituting the definition of $G$ and unfolding gives
$$
\Lambda_{\mathrm{U}}(G)(z)=\int_{\Gamma_0(7)\backslash\mathbb H}\Theta_{L,0}^{\mathrm{KM}}(\tau,z)\,\overline{f_7(\tau)}\,v^6\,\frac{du\,dv}{v^2}.
$$
Absolute convergence follows from the decay of $f_7$ at every cusp and the polynomial growth of the theta function.

Since $|A_L|=7^6$ and the signature factor is $1$, the $S$-transformation law gives
$$
(\Theta_{L,0}^{\mathrm{KM}}|_6W_7)(\tau,z)=7^3(\Theta_{L,0}^{\mathrm{KM}}|_6S)(7\tau,z)=\sum_{\lambda\in L^\vee}\varphi_{\mathrm{KM}}(\lambda,7\tau,z).
$$
Since scalar multiplication by $-i$ is unitary and fixes $\mathcal D_{V}$ pointwise, we have $\varphi_{\mathrm{KM}}(-ix,z)=\varphi_{\mathrm{KM}}(x,z)$.
For $x\in L$, the defining formula for $\varphi_{\mathrm{KM}}(x,\tau,z)$ therefore yields
$$
\varphi_{\mathrm{KM}}(\delta_K^{-1}x,7\tau,z)=e\bigl(uQ(x)\bigr)\varphi_{\mathrm{KM}}(-i\sqrt v\,x,z)=\varphi_{\mathrm{KM}}(x,\tau,z).
$$
Using $L^\vee=\delta_K^{-1}L$ in the preceding sum, we obtain $\Theta_{L,0}^{\mathrm{KM}}|_6W_7=\Theta_{L,0}^{\mathrm{KM}}$.

Finally, $W_7$ normalizes $\Gamma_0(7)$, so we may substitute $\tau'=W_7\tau$ in the unfolded integral.
Together with $f_7|_6W_7=-f_7$, the identity just proved yields
$$
\Lambda_{\mathrm{U}}(G)(z)=\int_{\Gamma_0(7)\backslash\mathbb H}(\Theta_{L,0}^{\mathrm{KM}}|_6W_7)(\tau,z)\,\overline{(f_7|_6W_7)(\tau)}\,v^6\,\frac{du\,dv}{v^2}=-\Lambda_{\mathrm{U}}(G)(z),
$$
and hence $\Lambda_{\mathrm{U}}(G)=0$.
\end{proof}


\printbibliography

\end{document}